\documentclass[preprint,review,12pt]{elsarticle}

\usepackage{graphicx}
 \usepackage{color}
 \usepackage{xcolor}
 \usepackage{cancel}

\usepackage{booktabs}

\usepackage[T1]{fontenc}
\usepackage{amsmath}
\usepackage{amstext,amsbsy,amsfonts}
\usepackage{graphicx,amssymb}

\usepackage[hang,small]{caption}
\usepackage{tikz}
\usepackage{pgfplots}

\usepackage{subfig}
\usepackage{subfloat}
\usetikzlibrary{pgfplots.groupplots,external}

\usepackage{multirow}
\usepackage{xspace}
\usepackage{enumitem}

\biboptions{comma,round,sort&compress}
\definecolor{blu}{RGB}{7, 40, 129}
\definecolor{rosso}{RGB}{168, 0, 0}
\definecolor{celeste}{RGB}{17, 170, 242}  
\definecolor{senape}{RGB}{236, 195, 33}
\definecolor{verde}{RGB}{0, 153, 76}

\pgfplotsset{h/.style={solid, black,thin}}
\pgfplotsset{h3/.style={dashed, black,thin}}
\pgfplotsset{r1/.style={solid, blu,thick}}
\pgfplotsset{d1/.style={r1, dashed}}
\pgfplotsset{r2/.style={r1,  verde}}
\pgfplotsset{d2/.style={r2, dashed}}
\pgfplotsset{r3/.style={solid, brown, thick}}

\pgfplotsset{C1IGA1/.style={brown, dashed, solid, thick}}
\pgfplotsset{C1IGA2/.style={C1IGA1, dashed}}

\pgfplotsset{C2IGA1/.style={rosso, dashed, solid, thick}}
\pgfplotsset{C2IGA2/.style={C2IGA1, dashed}}

\pgfplotsset{iso1/.style={rosso, solid, thin, mark= o,mark options={solid},mark size=1}}
\pgfplotsset{iso2/.style={iso1, dashed, mark= diamond,mark options={solid},mark size=1}}
\pgfplotsset{iper1/.style={blu, solid, mark=pentagon,mark options={solid},mark size=1}}
\pgfplotsset{iper2/.style={iper1, dashed, mark= triangle,mark options={solid},mark size=1}}

\pgfplotsset{Sconv/.style={xlabel=$log(h)$, ylabel=$log(E_\sigma)$, grid style={dotted, gray},grid = both,width=0.26\linewidth, height = 0.45\linewidth,
		xmin=-1, xmax=0.5,
		ymax=0,ymin=-7.5,
		legend style={at={(1.05,0.95)},font=\footnotesize,anchor=north west }},
	every axis y label/.style={at={(ticklabel* cs:1.1)}},
	every axis/.append style={font=\footnotesize}}

\pgfplotsset{Uconv/.style={xlabel=$log(\sqrt{dof})$, ylabel=$log(1-u/u_{ref})$, grid style={dotted, gray},grid = both,width=0.26\linewidth, height = 0.45\linewidth,
		legend style={at={(1.05,0.95)},font=\footnotesize,anchor=north west }},
	every axis y label/.style={at={(ticklabel* cs:1.1)}},
	every axis/.append style={font=\footnotesize}}

\pgfplotsset{SIGconv/.style={xlabel=$log(\sqrt{dof})$, ylabel=$log(1-\sigma/\sigma_{ref})$, grid style={dotted, gray},grid = both,width=0.26\linewidth, height = 0.45\linewidth,
		legend style={at={(1.05,0.95)},font=\footnotesize,anchor=north west }},
	every axis y label/.style={at={(ticklabel* cs:1.1)}},
	every axis/.append style={font=\footnotesize}}

\usetikzlibrary{calc}

\usetikzlibrary{pgfplots.colormaps} 
\pgfplotsset{colormap={parula}{
rgb=(0.2422,0.1504,0.6603)
rgb=(0.25039,0.165,0.70761)
rgb=(0.25777,0.18178,0.75114)
rgb=(0.26473,0.19776,0.79521)
rgb=(0.27065,0.21468,0.83637)
rgb=(0.27511,0.23424,0.87099)
rgb=(0.2783,0.25587,0.89907)
rgb=(0.28033,0.27823,0.9221)
rgb=(0.28134,0.3006,0.94138)
rgb=(0.28101,0.32276,0.95789)
rgb=(0.27947,0.34467,0.97168)
rgb=(0.27597,0.36668,0.9829)
rgb=(0.26991,0.3892,0.9906)
rgb=(0.26024,0.41233,0.99516)
rgb=(0.24403,0.43583,0.99883)
rgb=(0.22064,0.46026,0.99729)
rgb=(0.19633,0.48472,0.98915)
rgb=(0.1834,0.50737,0.9798)
rgb=(0.17864,0.52886,0.96816)
rgb=(0.17644,0.5499,0.95202)
rgb=(0.16874,0.57026,0.93587)
rgb=(0.154,0.5902,0.9218)
rgb=(0.14603,0.60912,0.90786)
rgb=(0.13802,0.62763,0.89729)
rgb=(0.12481,0.64593,0.88834)
rgb=(0.11125,0.6635,0.87631)
rgb=(0.09521,0.67983,0.85978)
rgb=(0.068871,0.69477,0.83936)
rgb=(0.029667,0.70817,0.81633)
rgb=(0.0035714,0.72027,0.7917)
rgb=(0.0066571,0.73121,0.76601)
rgb=(0.043329,0.7411,0.73941)
rgb=(0.096395,0.75,0.71204)
rgb=(0.14077,0.7584,0.68416)
rgb=(0.1717,0.76696,0.65544)
rgb=(0.19377,0.77577,0.6251)
rgb=(0.21609,0.7843,0.5923)
rgb=(0.24696,0.7918,0.55674)
rgb=(0.29061,0.79729,0.51883)
rgb=(0.34064,0.8008,0.47886)
rgb=(0.3909,0.80287,0.43545)
rgb=(0.44563,0.80242,0.39092)
rgb=(0.5044,0.7993,0.348)
rgb=(0.56156,0.79423,0.30448)
rgb=(0.6174,0.78762,0.26124)
rgb=(0.67199,0.77927,0.2227)
rgb=(0.7242,0.76984,0.19103)
rgb=(0.77383,0.7598,0.16461)
rgb=(0.82031,0.74981,0.15353)
rgb=(0.86343,0.7406,0.15963)
rgb=(0.90354,0.73303,0.17741)
rgb=(0.93926,0.72879,0.20996)
rgb=(0.97276,0.72977,0.23944)
rgb=(0.99565,0.74337,0.23715)
rgb=(0.99699,0.76586,0.21994)
rgb=(0.9952,0.78925,0.20276)
rgb=(0.9892,0.81357,0.18853)
rgb=(0.97863,0.83863,0.17656)
rgb=(0.96765,0.8639,0.16429)
rgb=(0.96101,0.88902,0.15368)
rgb=(0.95967,0.91346,0.14226)
rgb=(0.9628,0.93734,0.12651)
rgb=(0.96911,0.96063,0.10636)
rgb=(0.9769,0.9839,0.0805)
}}

\begin{document}

  \newcommand\undermat[2]{%
  \makebox[0pt][l]{$\smash{\underbrace{\phantom{%
    \begin{matrix}#2\end{matrix}}}_{\text{$#1$}}}$}#2}

\newcommand{\bhp}   {\mbox {\boldmath $\hat{p}$}}
\newcommand{\bhu}   {\mbox {\boldmath $\hat{u}$}}
\newcommand{\bdv}   {\mbox {\boldmath $\dot{v}$}}

\newcommand{\Div}           {\text{div}}
\newcommand{\tr}[1]         {\text{tr}\left( #1 \right)}
\newcommand{\gr}     {\mbox{\rm grad}}
\newcommand{\dv}     {\mbox{\rm div}}
\newcommand{\rot}    {\mbox{\rm rot}}
\newcommand{\T}      {\mbox{\rm T}}
\newcommand{\onehalf}  {\frac{1}{2}}
\newcommand{\asix}     {\frac{1}{6}}

\newcommand{\blam}  {\mbox {\boldmath $\lambda$}}
\newcommand{\bxi}   {\mbox {\boldmath $\xi$}}
\newcommand{\beps}  {\mbox {\boldmath $\varepsilon$}}
\newcommand{\bep}   {\mbox {\boldmath $\epsilon$}}
\newcommand{\bsig}  {\mbox {\boldmath $\sigma$}}
\newcommand{\bvphi} {\mbox {\boldmath $\varphi$}}
\newcommand{\bfi}   {\mbox {\boldmath $\varphi$}}
\newcommand{\bchi}  {\mbox {\boldmath $\chi$}}
\newcommand{\bdel}  {\mbox {\boldmath $\delta$}}
\newcommand{\bvrho}  {\mbox {\boldmath $\varrho$}}
\newcommand{\bTheta}{\mbox {\boldmath $\Theta$}}
\newcommand{\bPi}   {\mbox {\boldmath $\Pi$}}
\newcommand{\btheta}{\mbox {\boldmath $\theta$}}
\newcommand{\bpsi}  {\mbox {\boldmath $\psi$}}
\newcommand{\bdeps} {\mbox {\boldmath $\dot{\varepsilon}$}}
\newcommand{\bheps} {\mbox {\boldmath $\hat{\varepsilon}$}}
\newcommand{\bdsig} {\mbox {\boldmath $\dot{\sigma}$}}
\newcommand{\bhsig} {\mbox {\boldmath $\hat{\sigma}$}}
\newcommand{\balf}  {\mbox {\boldmath $\alpha$}}
\newcommand{\bfeta} {\mbox {\boldmath $\eta$}}
\newcommand{\invk}          {\frac {\scriptstyle 1}{\scriptstyle k}}
\newcommand{\bom}           {\mbox {\boldmath $\omega$}}
\newcommand{\bOm}           {\mbox {\boldmath $\Omega$}}
\newcommand{\bPhi}          {\boldsymbol{\Phi}}
\newcommand{\bPsi}          {\mbox {\boldmath $\Psi$}}
\newcommand{\bbet}          {\mbox {\boldmath $\beta$}}
\newcommand{\bkap}          {\mbox {\boldmath $\kappa$}}
\newcommand{\bfet}          {\mbox {\boldmath $\eta$}}
\newcommand{\bgam}          {\mbox {\boldmath $\gamma$}}
\newcommand{\bUpsilon}      {\mbox {\boldmath $\Upsilon$}}
\newcommand{\btau}          {\mbox {\boldmath $\tau$}}
\newcommand{\brho}          {\mbox {\boldmath $\varrho$}}
\newcommand{\balp}          {\mbox {\boldmath $\alpha$}}
\newcommand{\bphi}          {\mbox {\boldmath $ \phi$}}
\newcommand{\eps}           {\mbox {$\epsilon$}}
\newcommand{\bkt}           {\mbox {\boldmath $k_t$}}
\newcommand{\bzeta}          {\mbox {\boldmath $\zeta$}}

\newcommand{\bdu}   {\mbox {\boldmath $\dot{u}$}}
\newcommand{\bepsd} {\mbox {\boldmath $\tilde{\varepsilon$}}}
\newcommand{\bsigd} {\mbox {\boldmath $\tilde{\sigma$}}}
\newcommand{\epsd}  {\mbox {$\tilde{\varepsilon}$}}
\newcommand{\sigd}  {\mbox {$\tilde{\sigma}$}}
\newcommand{\stat}  {\mbox {\rm stazionario}}
\newcommand{\mSeq}  {\mathbb{E}_{eq}}
\newcommand{\mDcm}  {\mathbb{D}_{cm}}

\newcommand{\mE}    {\mathbb{E}}
\newcommand{\mP}    {\mathbb{P}}
\newcommand{\mS}    {\mathbb{S}}
\newcommand{\mK}    {\mathbb{K}}
\newcommand{\f}         {\mbox {\boldmath $f$}}
\newcommand{\bd}        {\mbox {\boldmath $d$}}
\newcommand{\bp}        {\mbox {\boldmath $p$}}
\newcommand{\bu}        {\mbox {\boldmath $u$}}
\newcommand{\bv}        {\mbox {\boldmath $v$}}
\newcommand{\br}        {\mbox {\boldmath $r$}}
\newcommand{\bs}        {\mbox {\boldmath $s$}}
\newcommand{\bK}        {\mbox {\boldmath $K$}}
\newcommand{\bx}        {\mbox {\boldmath $x$}}
\newcommand{\bz}        {\mbox {\boldmath $z$}}
\newcommand{\by}        {\mbox {\boldmath $y$}}
\newcommand{\bS}        {\mbox {\boldmath $S$}}
\newcommand{\bA}        {\mbox {\boldmath $A$}}
\newcommand{\bB}        {\mbox {\boldmath $B$}}
\newcommand{\bC}        {\mbox {\boldmath $C$}}
\newcommand{\bD}        {\mbox {\boldmath $D$}}
\newcommand{\bH}        {\mbox {\boldmath $H$}}
\newcommand{\bF}        {\mbox {\boldmath $F$}}
\newcommand{\bQ}        {\mbox {\boldmath $Q$}}
\newcommand{\bE}        {\mbox {\boldmath $E$}}
\newcommand{\bL}        {\mbox {\boldmath $L$}}
\newcommand{\bl}        {\mbox {\boldmath $\ell$}}
\newcommand{\bR}        {\mbox {\boldmath $R$}}
\newcommand{\bT}        {\mbox {\boldmath $T$}}
\newcommand{\bw}        {\mbox {\boldmath $w$}}
\newcommand{\bM}        {\mbox {\boldmath $M$}}
\newcommand{\bn}        {\mbox {\boldmath $n$}}
\newcommand{\bm}        {\mbox {\boldmath $m$}}
\newcommand{\bt}        {\mbox {\boldmath $t$}}
\newcommand{\bbe}       {\mbox {\boldmath $b$}}

\newcommand{\ba}            {\mbox {\boldmath $a$}}
\newcommand{\bb}            {\mbox {\boldmath $b$}}
\newcommand{\bc}            {\mbox {\boldmath $c$}}
\newcommand{\be}            {\mbox {\boldmath $e$}}
\newcommand{\bq}            {\mbox {\boldmath $q$}}
\newcommand{\bN}            {\mbox {\boldmath $N$}}
\newcommand{\bJ}            {\mbox {\boldmath $J$}}
\newcommand{\bg}            {\mbox {\boldmath $g$}}
\newcommand{\bG}            {\mbox {\boldmath $G$}}
\newcommand{\bV}            {\mbox {\boldmath $V$}}
\newcommand{\bW}            {\mbox {\boldmath $W$}}
\newcommand{\bI}            {\mbox {\boldmath $I$}}
\newcommand{\bi}            {\mbox {\boldmath $i$}}
\newcommand{\bj}            {\mbox {\boldmath $j$}}
\newcommand{\bk}            {\mbox {\boldmath $k$}}
\newcommand{\bU}            {\mbox {\boldmath $U$}}
\newcommand{\bP}            {\mbox {\boldmath $P$}}
\newcommand{\bX}            {\mbox {\boldmath $X$}}
\newcommand{\bY}            {\mbox {\boldmath $Y$}}
\newcommand{\bef}           {\mbox {\boldmath $f$}}
\newcommand{\btK}           {\mbox {\boldmath $\tilde{K}$}}

\newcommand{\campo}[1]      {{\Huge \textit{#1}}}
\newcommand{\dddotu}        {\stackrel{...}{u}}
\newcommand{\pu}            {\dot{u}}
\newcommand{\pBu}           {\dot{\bf u}}
\newcommand{\du}            {\delta u}
\newcommand{\veps}          {\varepsilon}
\newcommand{\Bw}            {{\bf m}}
\newcommand{\Bs}            {{\bf s}}
\newcommand{\Bwd}           {\dot{\bf m}}
\newcommand{\Bwt}           {\tilde{\bf m}}
\newcommand{\axial}[1]{\text{axial}\left[#1\right]}
\newcommand{\Rot}  [1]{\bR         \left[#1\right]}
\newcommand{\Log}  [1]{\text{log}  \left[#1\right]}
\newcommand{\spin} [1]{\text{spin} \left[#1\right]}
\newcommand{\Skew} [1]{\text{skew} \left[#1\right]}
\newcommand{\bzero} {\mbox {\bf 0}}
\newcommand{\bbeta} {\mbox {\boldmath $\beta$}}

\newcommand{\bvec}[2]{\begin{bmatrix}#1 \\ #2 \end{bmatrix}}

\newcommand{\Rp}[1]{\bR_{1}[\bvphi^{(k)}_{e#1}]}
\newcommand{\Rs}[2]{\bR_{2}[\bvphi^{(k)}_{e#1},\bvphi^{(k)}_{e#2}]}
\newcommand{\Rt}[3]{\bR_{3}[\bvphi^{(k)}_{e#1},\bvphi^{(k)}_{e#2},\bvphi^{(k)}_{e#3}]}
\newcommand{\Rq}[4]{\bR_{4}[\bvphi^{(k)}_{e#1},\bvphi^{(k)}_{e#2},\bvphi^{(k)}_{e#3},\bvphi^{(k)}_{e#4}]}

\newcommand{\Qp}[1]{\bQ_{1}[\bd_{e#1}]}
\newcommand{\Qs}[2]{\bQ_{2}[\bd_{e#1},\bd_{e#2}]}
\newcommand{\Qt}[3]{\bQ_{3}[\bd_{e#1},\bd_{e#2},\bd_{e#3}]}
\newcommand{\Qq}[4]{\bQ_{4}[\bd_{e#1},\bd_{e#2},\bd_{e#3},\bd_{e#4}]}

\newcommand{\revAdd}[1] { {\color{blue}{#1}} }
\newcommand{\revDel}[1] { {\color{red}{#1}} }

\newcommand{\com}  {\;,\;\;}

\newcommand{\fracs}[2]{\frac{#1}{#2}}

\newcommand{\Figref}[1] {Fig.~\ref{#1}}
\newcommand{\Tabref}[1] {Tab.~\ref{#1}}

\begin{frontmatter}
\title{Limit Analysis Approach for Optimal Reinforcement Design of Masonry Structures: the Weak Reinforcement Concept}

\author[poliba]{Aguinaldo Fraddosio  \corref{cor1}}
\author[poliba]{Celeste Lasorella}
\author[poliba]{Mario Daniele Piccioni}
\author[unina] {Elio Sacco}   

\cortext[cor1]{Corresponding author: aguinaldo.fraddosio@poliba.it  } 
\address[poliba]{Politecnico di Bari, dICAR, Via Edoardo Orabona, 4, 70125 Bari, Italy} 
\address[unina]{Universit\'a di Napoli Federico II, DiSt, Via Claudio 21, Napoli, Italy}

\begin{abstract}

In this paper, first the non-standard limit analysis 2D problem for reinforced masonry structures is suitably reformulated, giving a innovative computational framework for the evaluation of the strengthening effects of tensile reinforcements. This framework is used for designing optimized reinforcements for masonry structures aimed at improving the structural capacity, but without substantial changes in the typical behavior of masonry structural systems, since the reinforcement layout leaves the possibility of cracks opening. Therefore, both rocking motions and adaptation to support settlements remain allowed, and the transmission of too high forces to the unreinforced structural elements supporting the reinforced one is prevented. This way, a new strengthening strategy, called "weak reinforcement" approach is obtained. This strategy pursues the challenging goal of combining structural safety and conservation for architectural heritage: indeed, not only the historical masonry is preserved as much as possible, but also the original special "masonry-like" behavior. The design, in principle, of optimized strengthening intervention through the weak reinforcement approach is exemplified through the application to a trilithon under seismic loads and to an arch under a concentrated load. The capacity of the structure to adapt to support settlements is investigated through a suitable procedure for minimizing the total potential energy of the system.

\end{abstract}

\begin{keyword}
Masonry, Limit Analysis, Reinforcement, Optimal design, Historic Constructions, Foundations Settlements Adaptability.
\end{keyword}

\end{frontmatter}

\section{Introduction}

Masonry structures represent the architectural heritage of numerous countries and encapsulate centuries of historical, cultural, and technical knowledge. From ancient aqueducts and Romanesque churches to Renaissance palaces and early industrial buildings, masonry has served as the main construction material due to its widespread availability, durability, and high compressive strength. However, the structural behavior of masonry is inherently governed by its low tensile strength and brittle nature, which requires special consideration in both structural analysis and design \cite{Drysdale1999,Hendry2004,Como2016}.

Masonry structures have been conceived (or must be designed) to carry external loads — typically gravity loads — primarily through compressive stresses \cite{Heyman1997}; this is made possible by the geometric configuration of the structure. As such, these structures can be vulnerable to foundation settlements  \cite{Portioli2016,Tralli2020,Iannuzzo2025}, long-term deformations \cite{Ferretti2006,Cecchi2012,Sánchez-Beitia2017}, and seismic actions \cite{Tomazevic1999,Asteris2014,Lagomarsino2021}. The history test has shown that masonry structures can also exhibit a surprising capacity against the latter load conditions, unjustifiable through conventional structural analysis techniques, due to the ability to adapt through crack opening and activation of rocking behavior. However, beyond certain load intensities, this adaptation capacity may be insufficient and consequently collapse may occur. 

Since preserving masonry structures is of paramount importance, often strengthening interventions are unavoidable to ensure a suitable capacity to withstand both static and dynamic expected external loads. 
Anyway, the more recent guidelines on structural restoration of architectural heritage \cite{ICOMOS2024} not only insist on the principle of minimum intervention, but also highlight that one of the goal of the conservation is that of preserving the technical values of the construction, relative to the structural concepts, the construction techniques and the materials. Moreover, it is clearly stated that the strict application of seismic codes can lead to drastic and often unnecessary measures, and that in the strengthening of historical structures, it may be acceptable not to fully enforce the structural requirements conventionally prescribed by standards oriented to modern structures (concept of the structural improvement, opposed to the concept of structural retrofitting). These criteria have to be carefully respected for avoiding irreversible damage of the architectural heritage.

A broad spectrum of reinforcement techniques has been employed, ranging from time-tested traditional methods to more advanced and recent solutions. Conventional interventions, such as the use of steel ties, buttresses, and jack arches, remain the preferred choice, whether they are resolving, due to their compatibility with historical construction practices. Moreover, from the structural point of view, traditional strengthening methods respect the original structural behavior. In contrast, recent solutions are generally united by the concept that reinforcement is obtained by introducing tensile strength in the structures, since the low tensile strength is viewed as a weakness, rather than a peculiar feature, capable of revealing unexpected structural advantages if studied through suitable structural analysis approaches. In fact, the modern history of masonry strengthening interventions began with the addition of steel and/or reinforced concrete elements. More recently, advancements in materials science have led to the employ of externally bonded fiber-reinforced polymers (FRPs) and fiber-reinforced cementitious matrix (FRCM) composites, which offer high strength-to-weight ratios, improved durability, and ease of application; strengthening effects, bond properties and adhesion issues have been studied through theorical and experimental investigations \cite{Grande2008, Cancelliere2010, Carozzi2018, Alecci2016,Leone2017,Meriggi2022,Castellano2025}.

Despite the advantages recalled, the use of modern reinforcement techniques might have several drawbacks. Indeed, while locally effective, these interventions can introduce unintended consequences for the global structural behavior of the system. Specifically, they can create zones of excessive stiffness, thereby disrupting the continuity and integrity of the original structural response. In fact, reinforcement elements can introduce discontinuities in stiffness, which may hinder the development of natural deformation patterns or generate localized stress peaks. This issue is particularly critical in historic buildings, where preserving original load paths and failure mechanisms is often as essential as improving structural strength. Although reinforcements may delay or mitigate certain failure modes, they can also instigate more brittle or less favorable mechanisms. Moreover, substantial changes introduced in the overall structural behavior could obliterate the so important masonry capacity to adapt to settlements and to withstand seismic actions through rocking motions, at the end resulting in a reduction of the structural capacity and/or in unexpected damages due to unfavorable stress redistribution \cite{Cimellaro2011, Arcidiacono2016, Borri2019}. 
These drawbacks are substantially due to the fact that in the in the strengthened parts of the structure the reinforcements totally prevent blocks articulations, inducing a continuum body-like behavior, and that the counterpart of the increase of the structural capacity of the strengthened parts is the transmission of much higher forces to the unstrenghtened parts, like, e.g., it can happen when a strengthened masonry domes is connected to unstrenghtened pillars or support walls.

To address the above-described complexities introduced by the reinforcements in structural behavior, computational modeling has become an indispensable tool \cite{Grande2011, Grande2023, Castellano2023}. Anyway, describing masonry non-linear behavior is already quite complex, and the addition of reinforcements requires the description of non-linear interactions between masonry and reinforcements, and this is even more complex. Moreover, the structural models needed require the characterization of a large number of mechanical parameters that is generally experimentally challenging, if not infeasible. This makes very difficult, and substantially out of reach of professionals the application of proper analysis techniques for describing the structural behavior of reinforced masonry and for designing effective reinforcement interventions.

On the other hand, limit analysis has proven especially effective in assessing the ultimate load-bearing capacity and identifying potential failure mechanisms in unreinforced masonry structures. This approach does not require a detailed stress–strain characterization of the materials involved, but instead focuses on the limit state, making it particularly well-suited to systems with unilateral frictional contact, such as dry-joint or weakly bonded masonry. \cite{Block2006,Milani2006,Brandonisio2020,Fraddosio2020}. Limit analysis is effective even if, as is usually the case for historical constructions, the initial stress distribution in the structure is practically unknowable as the result of a complex story of events, and includes both static and kinematic formulations, allowing engineers to evaluate collapse loads and mechanisms efficiently. For the application of limit analysis to masonry constructions, along with continuum approaches, rigid block models and computational models have been developed \cite{Livesley1978, Alexakis2015, Chiozzi2017, Michiels2017, Fortunato2018, Iannuzzo2021, Nodargi2022, Nodargi2023}.

Recent developments in the application of limit analysis to masonry structures allows to include also conditions governed by non-associative evolution laws such as, e.g., frictional sliding at block joints. Indeed, while for associative behavior the extension of limit analysis to masonry is a well established matter \cite{Heyman1997, Delpiero1998}, in the non-associative case the application of standard procedure can lead to unconservative results in terms of capacity; moreover, non-uniqueness issues for the collapse load multiplier arise \cite{Sacchi1968,Radenkovic1961, Salençon1973, DeSaxcé1998}. Several authors proposed suitable approaches for tackling the limit analysis problem for masonry blocky structures governed by non-associative interfaces, capable of representing the actual role of friction in the collapse 
\cite{Gilbert2006, Trentadue2013, Portioli2014, Portioli2015, Portioli2016, Nodargi2019AVF, Mousavian2020, Hua2022, Cocchetti2024, Gagliardo2024, Rios2025, trentadue2025}. Since the non-standard problem is is more complicated from a mathematical programming standpoint, the problem is in some cases transformed in a suitable sequence of associative problems; moreover, the employ of effective optimization approach has been studied.

Moreover, the effects of seismic loads can be included in limit analysis by considering suitable static loads representative of the maximum inertia loads induced by the deisgn earthquake \cite{Nodargi2021}. 

Finally, also the effect of masonry reinforcements can be taken into account in the limit analysis by defining suitable admissible domains for the reinforced material, also in terms of space hosting admissible thrust lines or surfaces \cite{Caporale2006, Grande2008, Fabbrocino2015, Chiozzi2017, Nela2022}.
    
In this paper, by suitably reformulating the non-standard limit analysis problem for reinforced masonry structures, a innovative computational framework for the assessment of the strengthening effects of reinforcements consisting of the introduction of tensile strength at masonry intrados and/or at the extrados is presented. In addition, a strategy for an optimal design of such kind of reinforcements is proposed through quantitative evaluations allowing for determining the placement and quantity of reinforcements. Here, the "optimal design" of the reinforcement is intended as a special design format aimed at enhancing the structural performance up the desired level, but without compromising the adaptive behavior of masonry structural systems, since the reinforcement layout leaves the possibility of cracks opening. Moreover, such kind of reinforcement can avoid unexpected and too high force transmission to the unreinforced structural elements supporting the reinforced one. Reinforcement leaves some rigid-body degrees of freedom to the reinforced structures, rather acting in terms of introduction of suitable internal forces rather than of internal constraints. The above strategy represents a new contribution in the field, paving the way for radical innovations in the way masonry reinforcement is conceived, toward the challenging goal of combining structural safety and conservation for architectural heritage, and foreshadowing interesting theoretical, experimental and technological developments.
In short, for the above the proposed reinforcement strategy can be named "weak reinforcement" approach since the reinforcement design. 
Here, technological aspects for practically realizing it are not studied, but current reinforcement systems like, e.g., FRP and FRCM, can be reasonably arranged for achieving the weak reinforcement conditions, for example by adopting mechanical anchors.
The weak reinforcement appears in line with the criteria in \cite{ICOMOS2024}, minimizing the strengthening intervention, allowing for graduating the reinforcement for obtaining the desired structural capacity (that, for the structural improvement concept may be lesser than a new structure), and above all preserving as best as possible the original structural behavior.

In this paper, a comprehensive methodology for the optimal "weak" reinforcement design for two-dimensional masonry structures using non-standard limit analysis is presented. The structural model consists of rigid blocks connected by zero-thickness interfaces governed by frictional contact laws. Reinforcements are applied at the edges of the interfaces and are assumed to act in tension only, thereby accurately representing the uniaxial behavior of most FRCM and FRP systems. 
The proposed approach integrates both static and kinematic limit analysis theorems. 
In particular, a suitable format of Static Limit Analysis (SLA) aimed at maximizing the collapse load multiplier while minimizing reinforcement forces is proposed, while Kinematic Limit Analysis (KLA) simulates collapse mechanisms and can incorporate prescribed settlements to evaluate the adaptive capacity of the structure.

Moreover, an energy-based formulation based on the principle of total potential energy minimization is addressed for evaluating the effect of differential settlements on the structure, also in the presence of reinforcements.

By combining numerical limit analysis optimization approaches with mechanical modeling, this framework enables the identification of reinforcement configurations that not only improve structural capacity but also preserve the inherent mechanism formation capacity of masonry structures. This is particularly important for structures subjected to foundation settlements or differential displacements. Such characteristics are essential to preserve the isostatic adaptability of masonry systems, even in the presence of structural reinforcement intervention.
Indeed, classical reinforcement approaches based on the application of tensile resisting materials on masonry intrados and/or extrados, preventing hinges opening, in case of foundation settlements might lead to the introduction of high stress states capable of determining local or even local strength crisis of the material, totally canceling the adapting capacity to settlements typical of masonry structures.

The paper is summarized as follows. 
Section 2 is devoted to the position of the problem, and in Section 3 kinematic (Section 3.1) and equilibrium (Section 3.2) conditions are formulated for reinforced 2D blocky masonry structures. In Section 4, after presenting the general setting of the limit analysis problem for non-associative masonry block interfaces, two approaches are proposed for SLA (Section 4.1): the first is more conventional, whereas the second is more convenient for the purposes of optimal design of reinforcements. In Section 4.2, the procedure for KLA is presented. Finally, in Section 5 a total potential energy minimization approach is proposed to analyze masonry mechanisms under support settlements. 
 
In Section 6 the innovative weak reinforcement technique is outlined, and in Section 7 the proposed methodology is applied to two simple yet illustrative structural systems to demonstrate its effectiveness.
The first case (Section 7.1) examines a trilithon structure under live loads representative of maximum seismic inertia load in the presence of an incremental reinforcement quantity. Numerical results indicate a marked increase in the collapse load multiplier as the reinforcement increases, accompanied by a variation in failure mechanisms. Internal forces within each reinforcement element are systematically evaluated, enabling data-driven material optimization to enhance structural performance.
The second case (Section 7.2) involves a masonry arch under a vertical eccentric live load, for which experimental laboratory data are available for benchmarking. The arch is analyzed either for extrados or intrados reinforcement. The proposed analytical method demonstrates high predictive accuracy, with computed collapse loads and failure modes showing close agreement with experimental observations. Notably, while intrados reinforcements yield higher increase in the collapse load multiplier, they exhibit greater susceptibility to debonding under practical conditions. The optimization framework incorporates additional constraints to account for potential foundation settlements, ensuring robust reinforcement design.

\section{Position of the problem}

The design problem of "optimal" reinforcements for two-dimensional masonry structures is approached. Here, optimal means a special distribution of the reinforcement capable of significantly increasing the structural capacity of the reinforced part of the construction, but leaving the possibility of forming mechanisms for the latter, thus not violating the inherent nature of masonry structural behavior, and at the same time strongly reducing the actions transmitted to the unreinforced part of the construction.
The analysis is performed by introducing two basic assumptions regarding a) the masonry modeling and b) the type of analysis to be performed.

Masonry can be modeled by considering the real assembly of blocks joined by interfaces and eventually reinforcements. This kind of modeling approach seems to be particularly suitable for masonry, characterized by relatively large blocks with good mechanical properties, such as stiffness and strength. 
Blocks are joined by interfaces that may be completely dry, i.e., there is direct contact between one block and the adjacent one, or formed by thin layers of mortar. 
In the first case, the tensile strength of the interface is null, and the shear strength depends on the frictional properties of the contacting brick's faces, without cohesion; in the latter case, the tensile and shear strengths may vary greatly but anyway assuming very small values, which could also be difficult to derive from in situ experimental tests.
Therefore, a reasonable assumption in masonry modeling, mainly for very old constructions, is that both tensile strength and cohesion are negligible, i.e., zero. 
In instances where a mortar layer is present between blocks, and given that its thickness is considerably less than the dimensions of the blocks, it is generally accepted practice to assume a zero thickness for interfaces. This is achieved by appropriately expanding the size of the blocks until the mid-surface of the mortar joint is reached.
Consequently, both the dry and mortar interfaces between blocks can be modeled by assuming a unilateral contact law with friction. This approach has been widely used to investigate the non-linear mechanical response of masonry structures \cite{Sacco10}.

Notice that in cases like, e.g., brick masonry, for simplifying the analysis by reducing the number of interfaces, the blocks can represent groups of brick thought as macro-blocks.

Regarding the type of analysis to be performed to derive the information needed to design optimal reinforcement, capable of suitably improving the load-bearing capacity of a masonry structure, here the limit analysis is preferred to the evolutionary time-step approach because it is generally simpler, does not require the introduction of evolutionary constitutive laws, often very difficult to be experimentally characterized, and also because it can be considered particularly suitable for determining the structural capacity in problems governed by unilateral contact law with friction.

Therefore, in the following the limit analysis is adopted to derive the collapse load and the associated mechanism of masonry structures modeled by blocks joined by unilateral frictional interfaces.

The masonry structure is modeled by introducing $N_b$ blocks connected by $N_{\Im}$ interfaces and $N_R$ reinforcements.
As announced above, the analysis is developed in two-dimensional space, introducing a Cartesian coordinate system $(O,x,y)$ and employing a vectorial notation.

The blocks constituting the structure are characterized by arbitrary polygonal shapes, identified by the coordinates of its vertexes.
Interfaces are rectilinear, so the interface connecting two typical blocks, $b_j$ and $b_k$, is identified by four points, two belonging to $b_j$ and two to $b_k$. Since interfaces have zero thickness, the two pairs of points lying on the two boundaries of the blocks have exactly the same positions.
At the typical $i$-th interface $\Im_i$, the tangential unit vector $\bt^i$ and the normal unit vector $\bn^i$ are introduced, along with an abscissa $\xi$ with origin in the middle point of the interface that progresses in the positive tangential direction.

It is assumed that reinforcements can only be placed at the extremities of the each interface, so that $N_R\leq2\ N_{\Im}$, and that they act only along the direction perpendicular to the interface to whom belonging the point where the reinforcement is. It is evident that the latter option can be readily eliminated, provided that the direction of the reinforcement action is established.
This choice was made because reinforcements are generally not so effective for limiting shear sliding at interfaces, while they might work very well to limit the interface opening.
To clarify the notation used in the following, $r_i^{(1)}$ and $r_i^{(2)}$ identify the reinforcements associated with the $i$-th interface at the $\xi=-L_i/2$ and $\xi=L_i/2$ extremities, respectively. 
In this context, $L_i$ represents the length of the $i$-th interface.
It should be noted that the interface may be either unreinforced or reinforced at one or both of the extremities.
The force exerted by the reinforcement at acting in $r_i^{(p)}$, i.e. at the $p$-th extremity ($p=1,2$) of the $i$-th interface, is denoted as $R_i^{(p)}$.

The kinematics of the structure is defined by the displacement vector of the centroids of the blocks, i.e. 
$\bU=\left\{ \bU_1^{\T}\;\;  \bU_2^{\T}\;\; ...\;\; \bU_{N_b}^{\T} \right\}^{\T}$ with $\bU_j=\left\{U_j\;\; V_j\;\; \Phi_j\right\}^{\T} $ the vector containing the three degrees of freedom of the $j$-th block: 
translation along $x$, translation along $y$, rotation around the centroid.

Two classes of forces acting on the centroid of the blocks are considered:
\begin{itemize}
\item
the external forces, i.e.
\begin{itemize}
\item
the dead (or fixed) active forces $\bF^d$ ;
\item
the live active forces $\bF^{\ell}$ , i.e., forces variable according to a load multiplier $\lambda$,
\item
the reactive forces $\bF^V$ of the external constraints;
\end{itemize} 
\item
the internal forces applied at the centroid of the blocks, i.e.
\begin{itemize}
\item
the forces $\bF^{X}$ arising from the interfaces interactions  $\bX$,
\item
the forces transmitted by the reinforcements, denoted by $\bF^R$, due to the presence of individual reinforcements, each one exerting a force $\bR$.
\end{itemize}

\end{itemize}

All the forces introduced above are represented by vectors having three components for each block, which, for duality with the kinematics, represent the two force components directed along $x$ and $y$ and the moment acting on each block.

\section{Kinematics and equilibrium of the system of joined blocks}\label{sec:kin_eq}

\begin{figure}[ht!]
\centering
\includegraphics[scale=0.6]{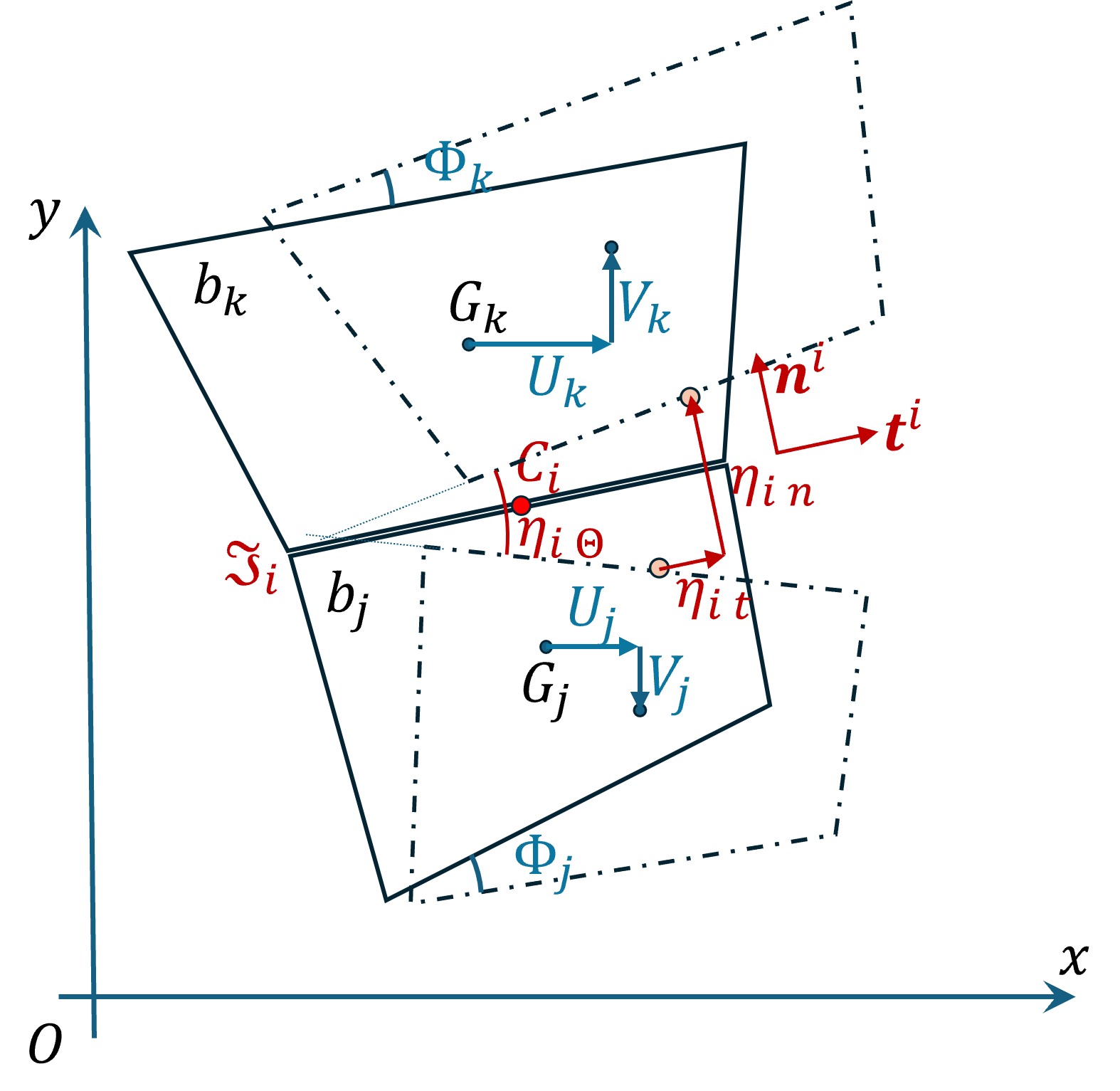}
\caption{Schematics of the kinematics of two blocks subjected to rigid motion: displacements of the blocks centroids and relative displacements at the interface. \label{fig:kin_2blocks}}
\end{figure}

In this section, the compatibility and equilibrium equations of the considered system of blocks joined by interfaces are stated aiming at formalizing the limit analysis (LA) problem.
Figure \ref{fig:kin_2blocks} illustrates the schematics of a possible rigid motion of two joined blocks.

\subsection{Kinematics}\label{sec:kin_eq_sub}

\begin{itemize}
\item
The kinematic compatibility equation relating the global displacement vector $\bU$ with the interface displacements is written for the whole structure as:
\begin{equation}
\bB\ \bU-\bfeta=\mathbf{0} \,,
\label{eq:kin_comp_int}
\end{equation}
where: \\
$\bfeta=\left\{ \bfeta_1^{\T}\;\;  \bfeta_2^{\T}\;\; ...\;\; \bfeta_{N_{\Im}}^{\T} \right\}^{\T}$ is the vector of the interface relative displacements with $\bfeta_i=\left\{\eta_{i\,t}\;\; \eta_{i\,n}\;\; \eta_{i\,\Theta}\right\}^{\T}$, whose components represent the interface relative displacement along $\bt^i$, along $\bn^i$ and the relative rotation, respectively;
\\
$\bB$ is the compatibility matrix able to express the interface relative displacement $\bfeta$ in terms of the block displacements. In particular, $\bB$ is the matrix obtained by assembling the matrices $\bB_i$, with $i=1,..,N_{\Im}$. For the typical interface $\Im_i$, joining the blocks $b_j$ and $b_k$, the compatibility matrix is given by:
\begin{equation}
\bB_i= \bQ_i \ \left[ \begin{array}{cccccc}
-1  &  0  &  y_i^C-y_j^G  &  1  &  0  &  -\left(y_i^C-y_k^G\right)\\
 0  &  \mathrm{-1} & -\left(x_i^C-x_j^G\right) & 0 & 1 & x_i^C-x_k^G\\
 0 & 0 & -1 & 0 & 0 & 1 \\
\end{array}
\right]
\label{eq:Bi} \,,
\end{equation}
with $\left(x_i^C,y_i^C\right)$ the coordinates of the interface center, $\left(x_j^G,y_j^G\right)$ and $\left(x_k^G,y_k^G\right)$ the coordinates of the centroids of $b_j$ and $b_k$, respectively, and  $\bQ_i$ the interface rotation matrix from the global  $(O,x,y)$ to the local coordinate system, referred to the interface directions $\bt^i$ and $\bn^i$.

The interface relative displacement is written as:
\begin{equation}
\bfeta = \bdel + \widetilde{\bN}\ \bbeta \,,
\label{eq:eta_betadelta}
\end{equation}
where $\bdel$ is an inelastic relative displacement (Volterra's distortion), prescribed at the interface, and 
$\bbeta$ is the contact vector satisfying the condition:
\begin{equation}
\bbeta \geq\mathbf{0} \,,
\label{eq:beta}
\end{equation} 
where by the inequality between vectors is here denoted in short the same inequality for each component of the vectors.

The vector $\bbeta$ is the assemblage of the vectors $\bbeta_i$ with $i=1,..,N_{\Im}$ interfaces, where the typical $i$-th vector, in components, is:
\begin{equation}
\bbeta_i=\left\{s_i\ \ {{\Delta u}_i}^+\ \ {{\Delta u}_i}^-\ \ {\Theta_i}^+\ \ {\Theta_i}^-\right\}^{\T}
\label{eq:betai},
\end{equation}
collecting the following kinematic parameters characterized by Eq. \eqref{eq:beta} to be always not negative: $s_i$ is the unilateral interface opening occurring along the normal direction $\bn^i$, ${\Delta u}_i$ is the interface sliding along the tangential direction  $\bt^i$, $\Theta_i$ is the relative rotation, while the superscripts $\cdot^+$ and $\cdot^-$ select the positive and negative part of the quantity to whom are referred, respectively, see Figure \ref{fig:beta}, where a schematic of the geometrical meaning of the components of the vector $\bbeta_i$ is illustrated.

\begin{figure}[ht!]
\centering
\includegraphics[scale=0.5]{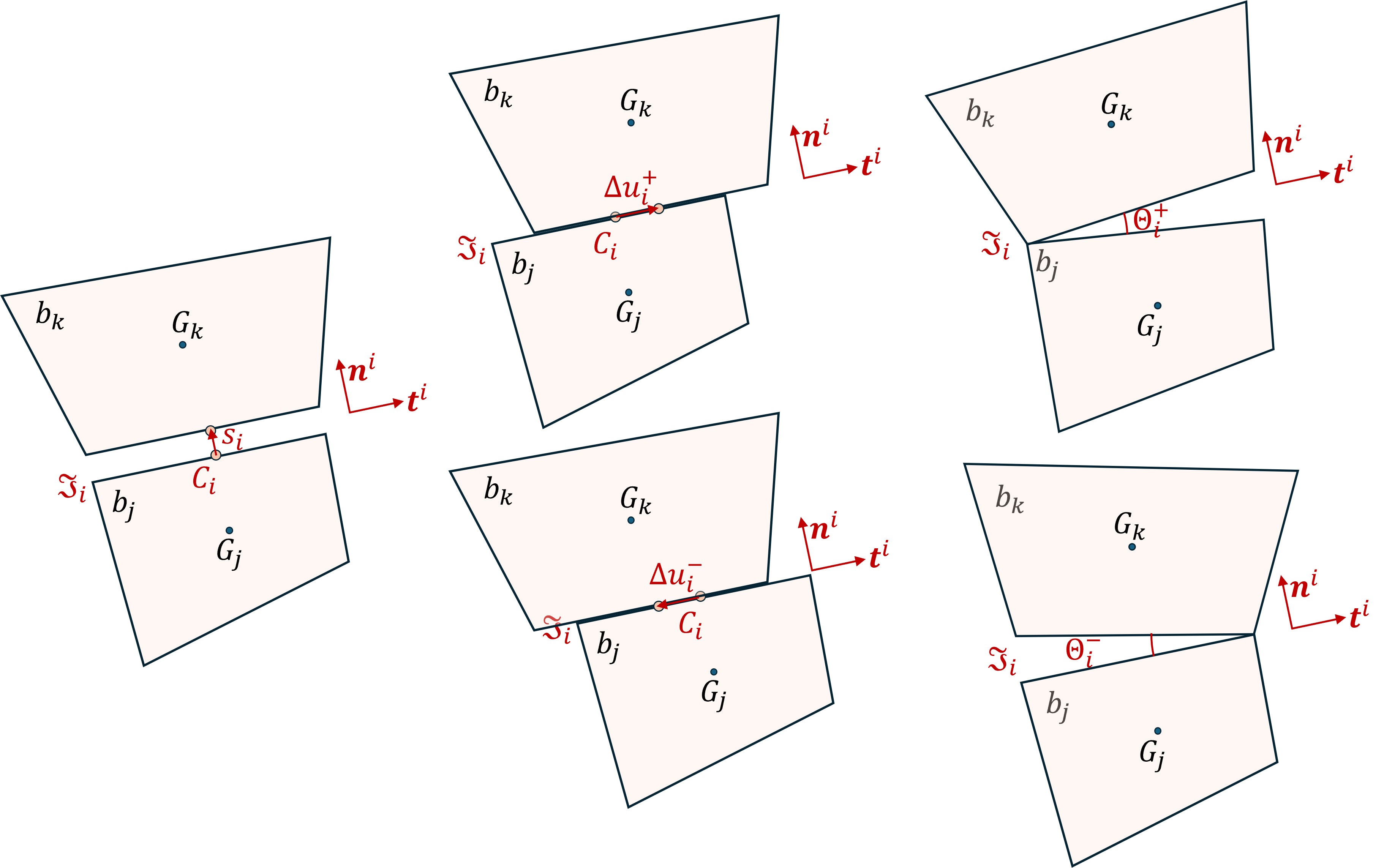}
\caption{Schematics illustrating the meaning of the components of  $\bbeta_i$. \label{fig:beta}}
\end{figure}

Finally, the matrix  $\widetilde{\bN}$, which governs the kinematic evolution law, is constructed by assembling the individual matrices $\widetilde{\bN}_i$, each corresponding to the $i$-th interface and defined as follows:
\begin{equation}
{\widetilde{\bN}}_i=\left[\begin{matrix}0&1&-1&0&0\\1&\widetilde{\mu}&\widetilde{\mu}&\rho_i&\rho_i\\0&0&0&1&-1\\\end{matrix}\right]
\label{eq:Nitilde},
\end{equation}
with $\widetilde{\mu}$ the dilatancy coefficient and $\rho_i$ half of the interface length, i.e. $\rho_i=L_i/2$. 

Note that, by following the approach indicated in \cite{sloan1995}, that turns out to be very useful for limit analysis of blocky structures \cite{Grillanda2025}, it si possible to write:
\begin{equation}
{\widetilde{\bN}}_i \ \bbeta_i=\left\{\begin{matrix}{{\Delta u}_i}^+-{{\Delta u}_i}^-\\s_i+\widetilde{\mu}\left({{\Delta u}_i}^++{{\Delta u}_i}^-\right)+\rho_i\left({\Theta_i}^++{\Theta_i}^-\right)\\{\Theta_i}^+-{\Theta_i}^-\\\end{matrix}\right\}=\left\{\begin{matrix}{\Delta u}_i\\{\Delta v}_i\\\Theta_i\\\end{matrix}\right\}
\label{eq:Nbetai}.
\end{equation}
with $\Delta u_i$ the contact tangent relative displacement along $\bt^i$, $\Delta v_i$ the contact normal relative displacement along $\bn^i$, and $\Theta_i$ the contact relative rotation, all components evaluated at the interface midpoint.

\item
The kinematic compatibility equation relating the global displacement vector $\bU$ with the displacement vector of the interface extremities $\bg$, activating the action of the reinforcements for the whole structure is written as:
\begin{equation}
\bC\ \bU = \bg \geq \mathbf{0}
\label{eq:kin_comp_T} \,,
\end{equation}
where: \\
$\bC$ is the compatibility matrix able to transform the block displacements into the relative displacements at the interface extremities, corresponding to the activation of the reinforcements. 
$\bC$ is a matrix obtained by assembling matrices that correspond to each single reinforcement. In particular, for the $r_i^{(p)}$-th reinforcement, associated with the $p$-th extremity (with $p=1,2$) of the $i$-th interface, connecting blocks $b_j$ and $b_k$, the local compatibility matrix can be written as:
\begin{equation}
\bC_i^{(p)} =
(\bn^i)^{\T}
\left[\begin{array}{rrrrrr}
-1 & 0 &  {y_i}^{(p)}-y_j & 1 & 0 & -{y_i}^{(p)}+y_k\\
 0 &-1 & -{x_i}^{(p)}+x_j & 0 & 1 &  {x_i}^{(p)}-x_k\\
\end{array}\right]
\label{eq:Ci} \,,
\end{equation}
with $\bn^i$ the normal vector to the $i$-th interface and $\left({x_i}^{(p)},{y_i}^{(p)}\right)$ the coordinates of the two nodes occupying the same position and joined by the reinforcement $r_i^{(p)}$.\\
Notice that the condition Eq. \eqref{eq:beta} implicitly implies that $\bg\geq\mathbf{0}$, as reported in Eq. \eqref{eq:kin_comp_T}.

\item
The kinematic compatibility equation relating the global displacement vector $\bU$ with the prescribed displacement vector $\bu$ of the constrained degrees of freedom is written as:
\begin{equation}
\bD\ \bU-\bu=\mathbf{0}
\label{eq:kin_constraint} \,,
\end{equation}
where: \\
$\bD$ is the matrix that selects the constrained displacements among all the displacements of the structure. This way it becomes easily possible to represent in the model eventual displacements due to imposed inelastic settlements of the structure's supports.

\end{itemize}

\subsection{Equilibrium}

For each of the blocks constituting the whole structure, the dead and live external forces, together with reactive forces and the interface and reinforcement internal forces, all applied at the centroid of the block, must be in equilibrium. This condition is expressed as:
\begin{equation}
\bF^{X}+\bF^R+\bF^V -\left(\bF^d + \lambda\ \bF^{\ell}\right) = \mathbf{0} 
\label{eq:equilibrium} \,.
\end{equation}

\begin{figure}[ht!]
\centering
\includegraphics[scale=0.6]{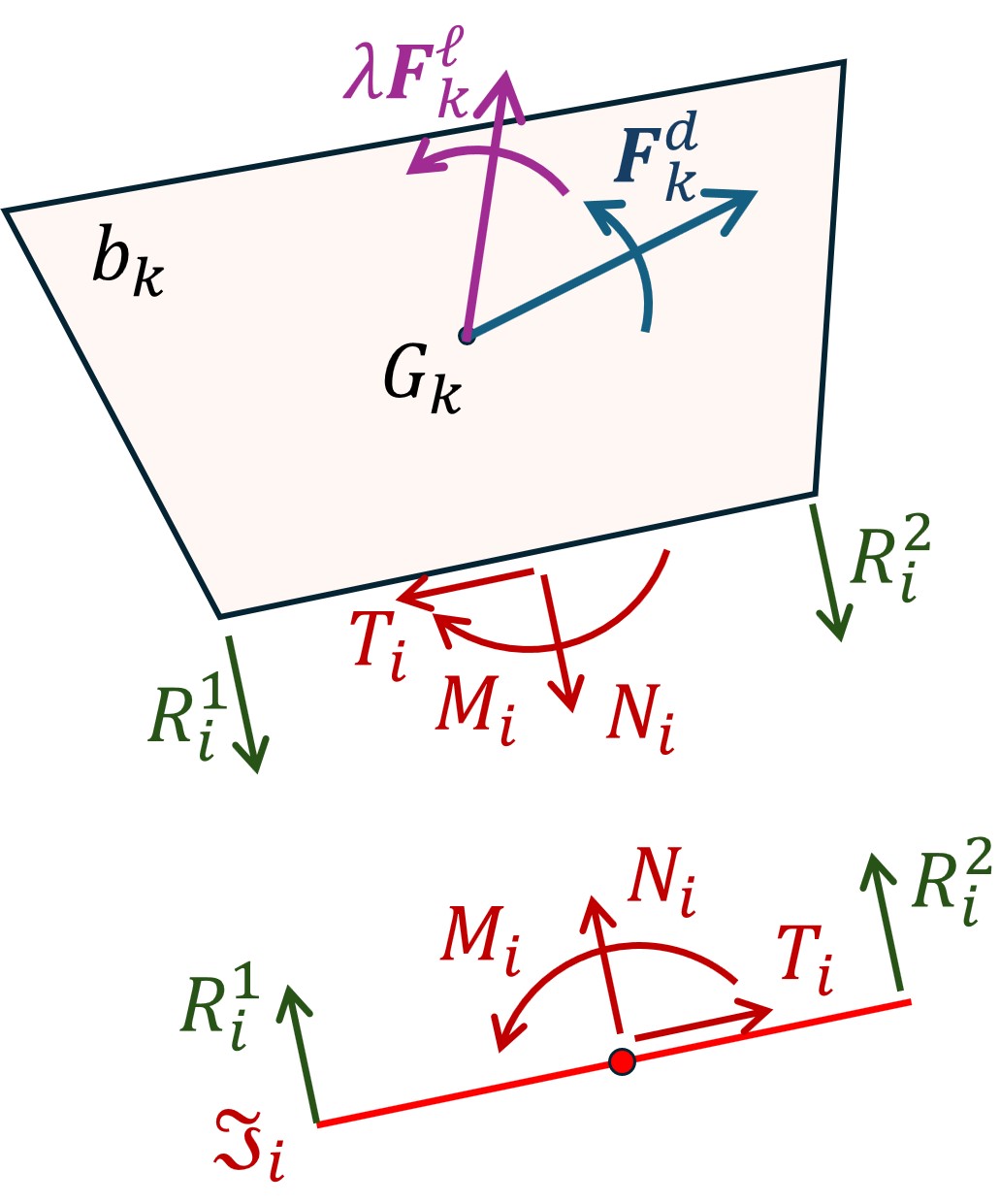}
\caption{Schematics illustrating the external forces, the interface tractions, and the tractions due to the presence of the reinforcements. \label{fig:forces}}
\end{figure}

In Figure \ref{fig:forces}, the external dead and live forces, the interface tractions, and the tractions resulting from the presence of the reinforcements are schematically represented. 
The subsequent section addresses the evaluation of the internal forces, i.e., the interface tractions, the reinforcement tractions, and reactive forces represented in Eq. \eqref{eq:equilibrium} by forces acting on the centroid of the generic block.

\begin{itemize}

\item
The interface forces $\bF^{X}$ applied at the centroid of the blocks are calculated as a function of the interface traction vector $\bX$, whose components are the resulting interaction forces and the moment, that each block transmits to the joined one through the interface. 
The typical traction vector evaluated at the $i$-th interface is thus $\bX_i=\left\{T_i \;\; N_i \;\; M_i\right\}^{\T}$, whose components represent the resulting shear force, normal force, and moment according to the local coordinate system defined by the interface tangent $\bt^i$ and normal $\bn^i$ unit vectors. 
By duality with the kinematic equation Eq. \eqref{eq:kin_comp_int}, the relationship between the interface forces $\bF^{X}$ and the interface traction vector $\bX$ is:
\begin{equation}
\bF^{X} = \bB^{\T} \ \bX
\label{eq:FBX}\,.
\end{equation}

Here, the interface traction vector has to satisfy the admissibility conditions: 

\begin{equation}
\bN^{\T} \ \bX \leq \mathbf{0}
\label{eq:NTXle0} \,,
\end{equation}
with $\bN$ a matrix, obtained by assembling matrices corresponding to each interface of the system. Eq. \eqref{eq:NTXle0} represents five inequalities enforcing at the interfaces the impossibility of tensile stresses, the compatibility of tangential forces with the friction available, and that the resulting normal force must lie within the interface. The typical matrix for the $i$-th interface takes the form:
\begin{equation}
\bN_i=\left[\begin{matrix}0&1&-1&0&0\\1&\mu&\mu&\rho_i&\rho_i\\0&0&0&1&-1\\
\end{matrix}\right]
\label{eq:Ni} \,,
\end{equation}
where $\mu$ is the friction coefficient.

\item
The reinforcement force vector $\bF^R$ is obtained as a function of the traction vector for the reinforcements $\bR$ by the relation:
\begin{equation}
\bF^R = \bC^{\T} \ \bR
\label{eq:FCT} \,,
\end{equation}
with the reinforcement force vector $\bR$ satisfying the condition:
\begin{equation}
\bR - \mathbf{1} \ R_y \leq \mathbf{0}
\label{eq:R} \,,
\end{equation}
being $\mathbf{1}$ a vector with all components equal to $1$ and $R_y$ the tensile yield limit force supported by the reinforcement.
Notice that by Eq. \eqref{eq:NTXle0}, Eq. \eqref{eq:R} also indicates that it is assumed that the reinforcements do not work in compression. 
The value of the tensile force in the reinforcement is equal to $R_y$ when it is in a limit state.

\item
The reactive force vector $\bF^{V}$ collects the vectors of the reactions of the constraints applied to each of the blocks, referred to the centroid of the blocks.
Therefore, by duality with the kinematic equation Eq. \eqref{eq:kin_constraint}, the reactive force vector associated with the kinematic constraints is given by the relationship: 
\begin{equation}
\bF^{V} = \bD^{\T} \ \bV
\label{eq:FBX2} \,.
\end{equation}

\end{itemize}

\section{Limit Analysis of Reinforced Masonry}\label{sec:LA}

Limit Analysis (LA) is a quite simple but powerful technique to evaluate the load-bearing capacity of a structure together with its corresponding collapse mechanism.
This last information can be very useful in the design of effective reinforcement systems: indeed, here it is used to introduce a proper reinforcement, the so-called "weak reinforcement", capable of inhibiting some mechanisms, allowing selectively the formation of new mechanisms characterized by a higher value of the collapse load.

At the end, the LA problem consists in solving the kinematic compatibility and equilibrium equations subject to the kinematic and static inequalities, expressing admissibility conditions, illustrated in the previous section, and linked by suitable complementary laws. 
In summary, in the general non-associative setting, the LA problem requires solving the following set of equations and inequalities, that can be conveniently combined to express the Lower Bound LA theorem (static approach) or the Upper Bound LA theorem (kinematic approach):
\begin{itemize}
\item
\textbf{kinematics}
\begin{equation}
\begin{array}{rrrr}
\qquad\qquad \qquad\qquad 
\bB\ \bU- \bdel - \widetilde{\bN}\ \bbeta &=& \mathbf{0}& \qquad\qquad \textrm{(a)} \\
\bD\ \bU-\bu      &=&    \mathbf{0} & \qquad\qquad \textrm{(b)} \\
\bbeta     &\geq & \mathbf{0} & \qquad\qquad \textrm{(c)} \\
\bC\ \bU         &\geq & \mathbf{0} & \qquad\qquad \textrm{(d)}
\end{array}
\label{eq:MCPkin} \,,
\end{equation}

\item
\textbf{equilibrium}
\begin{equation}
\begin{array}{rrrr}
\bB^{\T} \ \bX+\bD^{\T} \ \bV+ \bC^{\T} \ \bR -\left(\bF^d + \lambda\ \bF^{\ell}\right) &=&    \mathbf{0}& \qquad \textrm{(a)} \\
\bN^{\T} \ \bX    &\leq & \mathbf{0}  & \qquad \textrm{(b)} \\
\bR - \mathbf{1}\ R_y &\leq & \mathbf{0} & \qquad \textrm{(c)} 
\end{array}
\label{eq:MCPeq} \,,
\end{equation}

\item
\textbf{complementary laws}
\begin{equation}
\begin{array}{rrrr}
\qquad\qquad \qquad\qquad \qquad \qquad 
\bbeta^{\T} \; \bN^{\T} \ \bX &=& 0 & \qquad\qquad \textrm{(a)} \\
\left(\bC\ \bU \right)^{\T} \; \left(\bR - \mathbf{1}\ R_y \right) &=& 0 & \qquad\qquad \textrm{(b)}
\end{array}
\label{eq:MCPcomp} \,,
\end{equation}

\end{itemize}

where the kinematic equation equation Eq. \eqref{eq:MCPkin}$_{(a)}$ can be obtained by combining Eq. \eqref{eq:kin_comp_int} and Eq. \eqref{eq:eta_betadelta}, and the equilibrium equation Eq. \eqref{eq:MCPeq}$_{(a)}$ can be obtained by combining Eq. \eqref{eq:equilibrium}, Eq. \eqref{eq:FBX}, Eq. \eqref{eq:FCT}, and Eq. \eqref{eq:FBX2}.
The complementary equation Eq. \eqref{eq:MCPcomp}$_{(a)}$ ensures that if the contact vector of the $i$-th interface is such that $\beta_i> 0$, then the corresponding interface traction vector must be zero, i.e. $X_i=0$, and vice versa that
if the $i$-th interface traction vector is such that $X_i < 0$, then it must be $\beta_i=0$, i.e., no contact relative displacements at that interface occur.

Similarly, the complementary equation Eq. \eqref{eq:MCPcomp}$_{(b)}$ ensures that if the $i$-th component of the displacement vector of the interface extremities $\bg$ is $(\bC\ \bU)_i>0$, then the force in the reinforcements must reach the yielding limit, i.e. $R_i = R_y$, and vice versa that if at the $i-$th reinforcement
it results $R_i - R_y<0$, then no openings of that reinforcement are allowed, i.e., $(\bC\ \bU)_i=0$.

From a Mathematical Programming perspective, the problem governed by equations Eqs. \eqref{eq:MCPkin}-\eqref{eq:MCPcomp} represents a Nonlinear Mathematical Programming (NLP) problem, that is referred to as a Mixed Complementarity Problem (MCP). 
In particular, the considered problem is characterized by linear relationships, except for the nonlinear and nonconvex complementary conditions Eq. \eqref{eq:MCPcomp}$_{(a)}$, which make it difficult to solve.
Notice that the optimization problem governed by equations Eqs. \eqref{eq:MCPkin}-\eqref{eq:MCPcomp} becomes a Linear Programming (LP) problem when it is set $\widetilde{\mu}=\mu$, so that it results $\widetilde{\bN}=\bN$, i.e., when an associated frictional evolution law is considered.

Since the complementarity above discussed, the LA problem can be decoupled according to the static or the kinematic formulations, commonly used to solve the described optimization problem.

\subsection{Static Limit Analysis (SLA)}\label{sec:SLA}

The static formulation of the LA problem assumes that the unknowns are the force entities amplified by a load multiplier $\lambda$, and therefore it consists in finding the maximum value of the load multiplier $\lambda$ satisfying the equilibrium equations Eq. \eqref{eq:MCPeq} of the problem.
Two formulations of the static limit analysis are illustrated in the following. The first, called SLA1, considers the reinforcement force vector \bR\ as a further unknown, that consequently can be minimized having the objective, for example, of minimizing the quantity of reinforcing material. The second static limit analysis formulation, called SLA2, aims at minimizing the admissibility domain, and considers all the reinforcements working at their limit force $R_y$, as is usually done in static limit analysis of reinforced masonry.

\noindent
\textbf{SLA1} \\
The first formulation of the SLA can be formalized as the following optimization problem:
\begin{equation}
\begin{array}{llr}
\displaystyle{\max_{\lambda \;{ \bX \; \bV \; \bR}}}  & \displaystyle{\left\{\; \lambda  -
\frac{\alpha}{N_R\,R_y}\sum_{i=1}^{N_R}R_i \right\}} & \\
& \bB^{\T} \ \bX+\bD^{\T} \ \bV+ \bC^{\T} \ \bR -\left(\bF^d + \lambda\ \bF^{\ell}\right) = \mathbf{0} &\qquad (a) \\
& \bN^{\T} \ \bX \leq \mathbf{0} &\qquad (b) \\
& \bR - \mathbf{1}\ R_y \leq  \mathbf{0} &\qquad (c) 
\end{array}
\label{eq:stat1} \,.
\end{equation}
It is important to note that the maximization of the objective function contains not only the live load multiplier, denoted by $\lambda$, but also the value of the forces in the reinforcements normalized with respect to their tensile limit $R_y$ multiplied by the ratio $\alpha/(N_R R_y)$ with $\alpha$ a factor that weights the influence of the reinforcement forces in the maximization of the load multiplier. Notice that by setting a value for $\alpha$ in the interval $[0,\lambda]$ it is possible to optimize the reinforcement design, minimizing the quantity of reinforcing material. In particular, a suitable way for performing this optimized design is that of choosing for $\alpha$ a small value, suitably lesser than the value of $\lambda$ for the unreinforced structure.
The forces in the reinforcements are bounded below because of their unilateral response only in tension, and above because of the limit strength of the reinforcement material. This kind of limit analysis formulation allows for the derivation of the static collapse multiplier with minimal values of the forces in the reinforcements.
\\
The optimization problems Eq. \eqref{eq:stat1} can be conveniently transformed in a sequence of LP problems, properly modifying the inequality Eq. \eqref{eq:stat1}$_{b}$ as:
\begin{equation}
{\widetilde{\bN}}^{\T} \ \bX +\Delta\bN^{\T} \ \bX^{\star}  \leq \mathbf{0}  \label{eq:MCP15mod} \,,
\end{equation}
where the vector $\bX^{\star}$ represents the value of $\bX$ at the previous iteration and $\Delta\bN=\bN-\widetilde{\bN}$.
This sequential iterative approach is required in the general case of non-associative problems. Note that the related associated problem is simply solved as the first of the above iteration, by setting $\widetilde{\bN}=\bN$.
\\
\\
\noindent
\textbf{SLA2} \\
Suitable variations to the optimization problem Eq. \eqref{eq:stat1} can be implemented.
In fact, the presence of the reinforcements can be also accounted by modifying the interface admissibility law.
Simple equilibrium considerations allow for deriving the following equations for the typical $i$-th interface:
\begin{itemize}
\item
$r_i^{(1)}$ reinforcement:
\begin{equation}
\begin{array}{lll}
\;\;M_i+N_i \,\rho_i                &\leq & 0 \\
-M_i+N_i \,\rho_i -2 \rho_i\ R_y &\leq & 0
\end{array}
\label{eq:leftR} \,,
\end{equation}
\item
$r_i^{(2)}$ reinforcement:
\begin{equation}
\begin{array}{lll}
\;\;M_i+N_i \,\rho_i -2 \rho_i\ R_y &\leq & 0 \\
-M_i+N_i \,\rho_i                &\leq & 0
\end{array}
\label{eq:rightR} \,.
\end{equation}

\end{itemize}

Thus, conditions Eq. \eqref{eq:MCPeq}$_2$ can be rewritten as:
\begin{equation}
\bN^{\T} \ \bX - \widehat{\bR}  \leq \mathbf{0}
\label{eq:NTX_R} \,,
\end{equation}
where the vector $\widehat{\bR}$ is obtained by assembling the vectors $\widehat{\bR}_i$ with $i=1,2,..,N_{\Im}$, whose typical expression for the $i$-th interface is:
\begin{equation}
\widehat{\bR}_i = 2\rho_i\,R_y \,\left\{
\begin{array}{ccc}
\left\{0 \; \; \; 0 \; \; \; 0  \; \; \; 0 \; \; \; 1 \right\}^{\T} & \qquad \textrm{for reinforcement} & r_i^{(1)} \\
\left\{0 \; \; \; 0 \; \; \; 0  \; \;  \; 1 \; \; \; 0 \right\}^{\T} & \qquad \textrm{for reinforcement}& r_i^{(2)} 
\end{array} \right. \,.
\end{equation}
Thus, taking into account the above equations, the optimization problem Eq. \eqref{eq:stat1} of the LA for a 2D reinforced block masonry structure can be rewritten in the form:
\begin{equation}
\begin{array}{llr}
\displaystyle{\max_{\lambda \; \bX \; \bV}}  & \displaystyle{\left\{\; \lambda  \; \right\}} & \\
& \bB^{\T} \ \bX+\bD^{\T} \ \bV -\left(\bF^d + \lambda\ \bF^{\ell}\right) = \mathbf{0} &\qquad (a) \\
& \bN^{\T} \ \bX - \widehat{\bR} \leq \mathbf{0} &\qquad (b)  
\end{array}
\label{eq:stat2}\,,
\end{equation}
where the force vector $\bR$ does not appear explicitly because now the forces exerted by the reinforcements are included in the admissibility interface law.
\\
Even in this case, the nonlinear optimization problems Eq. \eqref{eq:stat2} and Eq. \eqref{eq:stat2} can be conveniently transformed in a sequence of LP problems, properly modifying the inequality Eq. \eqref{eq:stat2}$_{b}$ as:
\begin{equation}
{\widetilde{\bN}}^{\T} \ \bX +\Delta\bN^{\T} \ \bX^{\star} -\widehat{\bR}  \leq \mathbf{0} \label{eq:MCP25mod}\,,
\end{equation}
where the vector $\bX^{\star}$ represents the value of $\bX$ at the previous iteration and $\Delta\bN=\bN-\widetilde{\bN}$.
Note that the related associated problem is simply solved as the first of the above iteration, by setting $\widetilde{\bN}=\bN$.
\\
\noindent
\textbf{Remarks} \\
An interesting case can also be obtained by considering the yield strength $\ R_y$ of the interface as infinite. This assumption can be adopted for simplifying the problem when the reinforcement strength is much greater than other forces involved. 
In such a case, the inequality Eq. \eqref{eq:stat1}$_{(c)}$ disappears, as it is always satisfied, and taking into account the complementary law Eq. \eqref{eq:MCPcomp}$_{(b)}$, the compatibility equation Eq. \eqref{eq:MCPkin}$_{(d)}$ becomes: $\bC\ \bU=\mathbf{0}$.

\subsection{Kinematic Limit Analysis (KLA)}\label{sec:KLA}

Taking into account Eqs. \eqref{eq:MCPkin}-\eqref{eq:MCPcomp}, the kinematic formulation of the LA problem is derived by minimizing a suitable objective function derived from Eq. \eqref{eq:MCPkin}$_{(a)}$ by evaluating the work performed by the forces involved. 
In fact, premultiplying the equilibrium equation Eq. \eqref{eq:MCPeq}$_{(a)}$ by the displacement vector $\bU$, and assuming that the reinforcement force reaches its maximum value $R_y$, the optimization problem becomes:
\begin{eqnarray}
\min_{\bU \; \bbeta \; \ \bX \; \ \bV} & & \left\{
\bU^{\T}\left( \bB^{\T} \bX+ \bc\,R_y - \bF^d \right)+\bu^{\T} \bV  
\right\} \label{eq:kin} \\
&  & \quad  \bB\ \bU-\bdel-\widetilde{\bN}\ \bbeta = \mathbf{0}  \nonumber\\
&  & \quad  \bU^{\T}\ \bF^{\ell}-1 = 0  \nonumber \\
&  & \quad  \bbeta \geq \mathbf{0} \nonumber \\
&  & \quad  \bC\bU \geq \mathbf{0} \nonumber 
\label{eq:kin1},
\end{eqnarray}
where $\bc=\bC^{\T}\mathbf{1}$ and  Eq. \eqref{eq:kin1}$_{(3)}$ indicates that the function to be minimized is the load multiplier. 
This is a (complex) nonlinear optimization problem and, moreover, it appears as spurious because the unknowns vectors are not purely kinematic entities. In fact, the interface traction vector $\bX$ and the reactions vector $\bV$ are also unknowns of the problem.
Of course, the term involving the reactions vector $\bV$ disappears when the prescribed displacement vector at the constrained degrees of freedom is equal to zero, i.e., when $\bu=\mathbf{0}$ (perfect constraints).
Moreover, based on the considered type of reinforcements, capable only of exerting forces orthogonal to the interface and applied at the interface extremities, the condition $\bC\bU\geq\mathbf{0}$ results implicitly included in the other conditions, provided that it is evident that where reinforcements are applied also distortions $\bdel$ cannot be applied. 

As it has been done for the Static Limit Analysis, the nonlinear term due to the non-associativity of the interface friction law can be treated by transforming the problem into a sequence of linear programming optimization problems. In fact, it holds true:
\begin{equation}
\begin{split}
\bU^{\T} \ \bB^{\T}\bX &=
(\bB\ \bU)^{\T} \ \bX \\
&= \left(\bdel+\widetilde{\bN}\bbeta \right)^{\T} \ \bX \\
&= \left(\bdel+\bN\bbeta-\Delta\bN\bbeta \right)^{\T} \ \bX\\
&= \left(\bdel-\Delta\bN\bbeta \right)^{\T} \ \bX
\end{split},
\end{equation}
where Eqs. \eqref{eq:MCPkin}$_{(a)}$ and \eqref{eq:MCPcomp}$_{(a)}$ have been taken into account, and $\Delta\bN$ is the differences between $\bN$ and $\widetilde{\bN}$, accounting for the difference between the associative and non-associative kinematic evolution law. 
Therefore, the  iterative KLA problem  for perfect constraints conditions can be written in the form:
\begin{eqnarray}
\min_{\bU \; \bbeta} & & \left\{
\bdel^{\T}\bX^{\star}-(\bbeta^{\star})^{\T}\Delta\bN^{\T}\bX^{\star}+
\bU^{\T}\left(\bc\,R_y - \bF^d \right)  
\right\} \label{eq:kin2} \\
&  & \quad  \bB\ \bU-\bdel-\bN\ \bbeta+ \Delta\bN\ \bbeta^{\star}= \mathbf{0}  \nonumber\\
&  & \quad  \bU^{\T}\ \bF^{\ell}-1 = 0  \nonumber \\
&  & \quad  \bbeta \geq \mathbf{0} \nonumber,
\end{eqnarray}
where the interface kinematic parameter vector $\bbeta^{\star}$ and traction vector $\bX^{\star}$  are the values computed at the previous iteration. Note that $\bX^{\star}$ is determined as the Lagrange multiplier of the constraint Eq. \eqref{eq:MCPkin}$_{(a)}$.

\section{Total potential enenrgy approach}

It would be a subject of interest to determine the configuration of the system of rigid blocks joined by frictional interfaces, and eventually in presence of reinforcements, under a given loading condition and subjected to prescribed inelastic relative displacements at some interface. This would involve simulating supports displacements.
In this case, the live forces are determined to be equal to the product of an assigned value of the loading multiplier, denoted as $\lambda_a$, and the reference forces $F^{\ell}$, resulting as: $\lambda_a F^{\ell}$. In what follows it is mandatory that the value of $\lambda_a$ is less than the value of the collapse multiplier.

In the spirit of developing an iterative method for the purpose of addressing an associated type problem at each iteration step, it is possible to introduce the total potential energy at each step. In the specific case, this energy takes the following form:
\begin{equation}
\Pi(\bU) = -(\Delta\bN\ \bbeta^{\star})^{\T}\bX^{\star}+\bU^{\T}\bc\ R_y -\bU^{\T} \left(\bF^d+\lambda_a \bF^{\ell}\right) 
\label{eq:TPE}
\end{equation}
where $\bX^{\star}$ and $\bbeta^{\star}$ are assumed to be known, the reinforcement force reaches its maximum value $R_y$.
The equilibrated solution can be obtained by minimizing the total potential energy under admissibility kinematic conditions:
\begin{eqnarray}
\min_{\bU \; \bbeta} & &  \Pi(\bU) \label{eq:minTPE} \\
&  &  \quad \bB\ \bU-\bdel-\bN\ \bbeta+ \Delta\bN\ \bbeta^{\star}= \mathbf{0}  \nonumber\\
&  & \quad \bbeta \geq \mathbf{0} \nonumber,
\end{eqnarray}
Note that, at the typical iteration, the values of $\bX^{\star}$ and $\bbeta^{\star}$ are the ones computed at the previous iteration.

\section{The weak reinforcement optimal design}\label{sec:OR}

For the design of the weak reinforcement, the two approaches to the limit analysis problem, static and kinematic, can be suitably exploited. Indeed, each of these approaches have specific characteristics. 
In fact, the proposed static formulations, SLA1 and SLA2, allow for force minimization in the reinforcements, while the load multiplier $\lambda$ is maximized.
On the other hand, kinematic analysis allows for assessing the collapse mechanism and for the introduction of distortions within the interfaces, that can easily simulate the presence of settlements.
Therefore, the specific features of both LA approaches, together with analysis of the system of joined blocks through the minimization of the total potential energy, can be conveniently used to define a procedure for the optimized design of the masonry reinforcement.

A relevant aspect to be highlighted, is one of the principal advantages of masonry structures in comparison to concrete and steel structures: its capacity of articulating by joints openings, forming mechanisms under seismic actions, undergoing rocking motions without collapsing, or adapting to even significant foundation settlements. 
This enables the structure to achieve new dynamic and static equilibrium configurations without the occurrence of overstress.
The weak reinforcement concept is a strategy for introduction of reinforcing elements on a masonry structure aiming at reasonably increase the maximum load capacity of the structure without, however, affecting its remarkable ability to articulate under seismic actions and to regain different equilibrium configurations after foundation settlements, without inducing overstresses. Notice that the latter capabilities of masonry structures, that can be considered distinctive features of the masonry as a structural materials, and are related to the low tensile strength, are substantially lost with conventional reinforcement techniques based on the application of tensile resisting strengthening materials. 
Therefore, the weak reinforcement approach can preserve the original structural behavior of masonry, and this is consistent with conservation requirements. Furthermore, this approach is useful also for realizing another conservation requirement, i.e., that the reinforcement applied to the structure must be minimal, not greater than the quantity that is actually required for obtaining the structural capacity needed.

The proposed limit analysis formulations, in conjunction with the utilization of total potential energy, can serve as valuable tools for the optimal design of weak reinforcements for masonry structures.
In fact, the optimization of the reinforcement system can be carried out in three steps:
\begin{itemize}
\item  firstly, the SLA approach (in particular, SLA1)  can be employed to define a minimal weak reinforcement system;
\item secondly, this system can then be validated using the KLA approach, which involves inelastic relative displacements at the interfaces and thus can simulate foundation settlements;
\item thirdly, the structural integrity of the reinforced structure can be assessed through analysis of total potential energy minimization. This analysis enables the evaluation of the system's capacity to regain a new configuration following foundation settlements without inducing overstresses.
\end{itemize}

By adopting the proposed structural analysis approach, the introduction of reinforcing elements on the masonry structure may be undertaken in a stepwise manner, with each element being incorporated individually.
As previously stated, the considered reinforcements are not effective in limiting eventual shear sliding. However, they are highly effective in limiting interface openings on either side. Therefore, a strategy for determining the optimal location where adding a reinforcement is to evaluate where the maximum interface opening occurs, and then to introduce a reinforcement at that point.
This process is accomplished by introducing one reinforcement at a time and reanalyzing the structure: if the obtained load bearing capacity of the structure is still insufficient, the design process should be repeated by evaluating the maximum interface opening in presence of the previously introduced reinforcements, until the desired structural capacity is reached.

\section{Numerical results}
In the following, the weak reinforcement concept and the limit analysis format here developed for its design will be discussed in the light of application to simple structures, representative of a large class of applicative problems. Along with the improvements in the load bearing capacity, also the residual capacity of forming mechanism and of adapting to support settlements after the reinforcement will be investigated.

\subsection{The trilite}

\begin{figure}[ht!]
\centering
\includegraphics[scale=0.8]{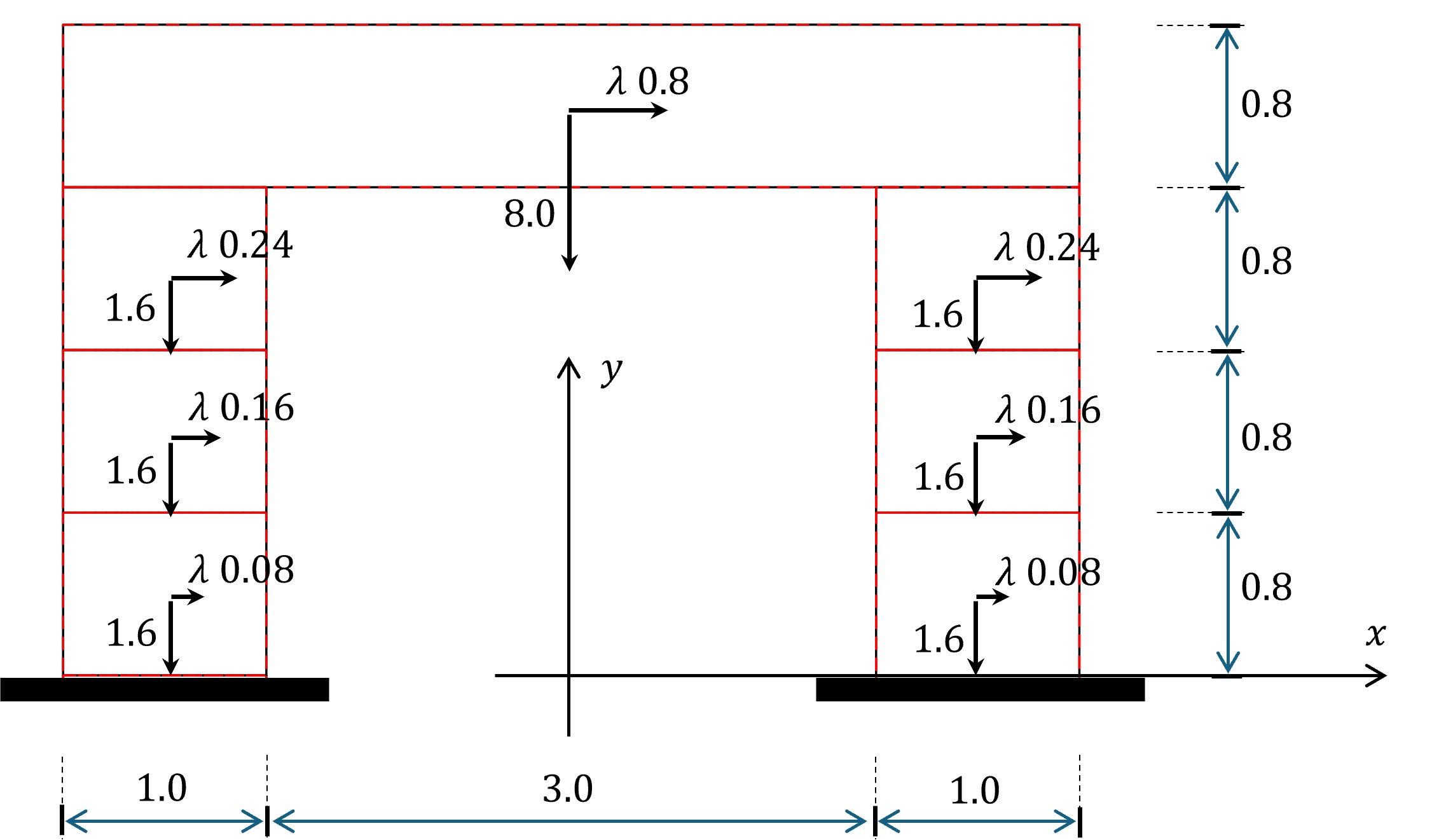}
\caption{Kinematics and collapse multiplier for the trilithon with increasing number of reinforcements. \label{fig:trilite}}
\end{figure}

The trilithon is a simple structural system that can serve as the first useful model for understanding the reinforcement optimization proposed herein. Specifically, a trilithon structure is here schematized as seven joined blocks, as illustrated in Figure \ref{fig:trilite}. The geometrical dimensions, loading conditions, and constraints of the model are also reported therein.
It should be noted that the application is intended to serve as a purely illustrative guide, to elucidate the proposed procedure. As such, it does not include any units of measurement.
The weights, and therefore the dead forces, are evaluated assuming a mass density $\gamma=2$ per unit area, and horizontal live forces representing a linear distribution of inertia forces increasing with the height are described by multiplying the mass of each block by a coefficient $\omega$ which varies linearly along the vertical direction according to the formula $\omega=(y+0.4)/8$. Clearly, these horizontal forces represents the maximum inertia forces on the structure induced by a seismic events.
To study the seismic capacity of the trilithon, these horizontal forces are considered as live forces, and therefore are multiplied by a multiplier $\lambda$.

The friction and the dilatancy coefficients at the interfaces are set as $\mu=1$ and $\tilde{\mu}=0$, respectively, while the strength of the reinforcement is assumed to be $R_y=10$. 
In the application of the SLA1, it is set $\alpha=0.1$, where (see Eq. \eqref{eq:stat1}) $\alpha$ weights the influence of the forces exerted by the reinforcements in the maximization of the collapse load multiplier; in particular, the value chosen for $\alpha$ is about $1/25$ of the collapse load multiplier for the unreinforced trilite; a greater value of $\alpha$ would have lead to a greater influence of the reinforcement on the structural capacity. 

The collapse load and the corresponding kinematics can be calculated by applying either the static theorem (both SLA1 and SLA2 approaches) or the kinematic theorem (KLA). In fact, all three possibilities were used, and the evaluation of the collapse load and the corresponding activated kinematics yielded the same result in all cases.

\begin{figure}[ht!]
\centering
\includegraphics[scale=0.8]{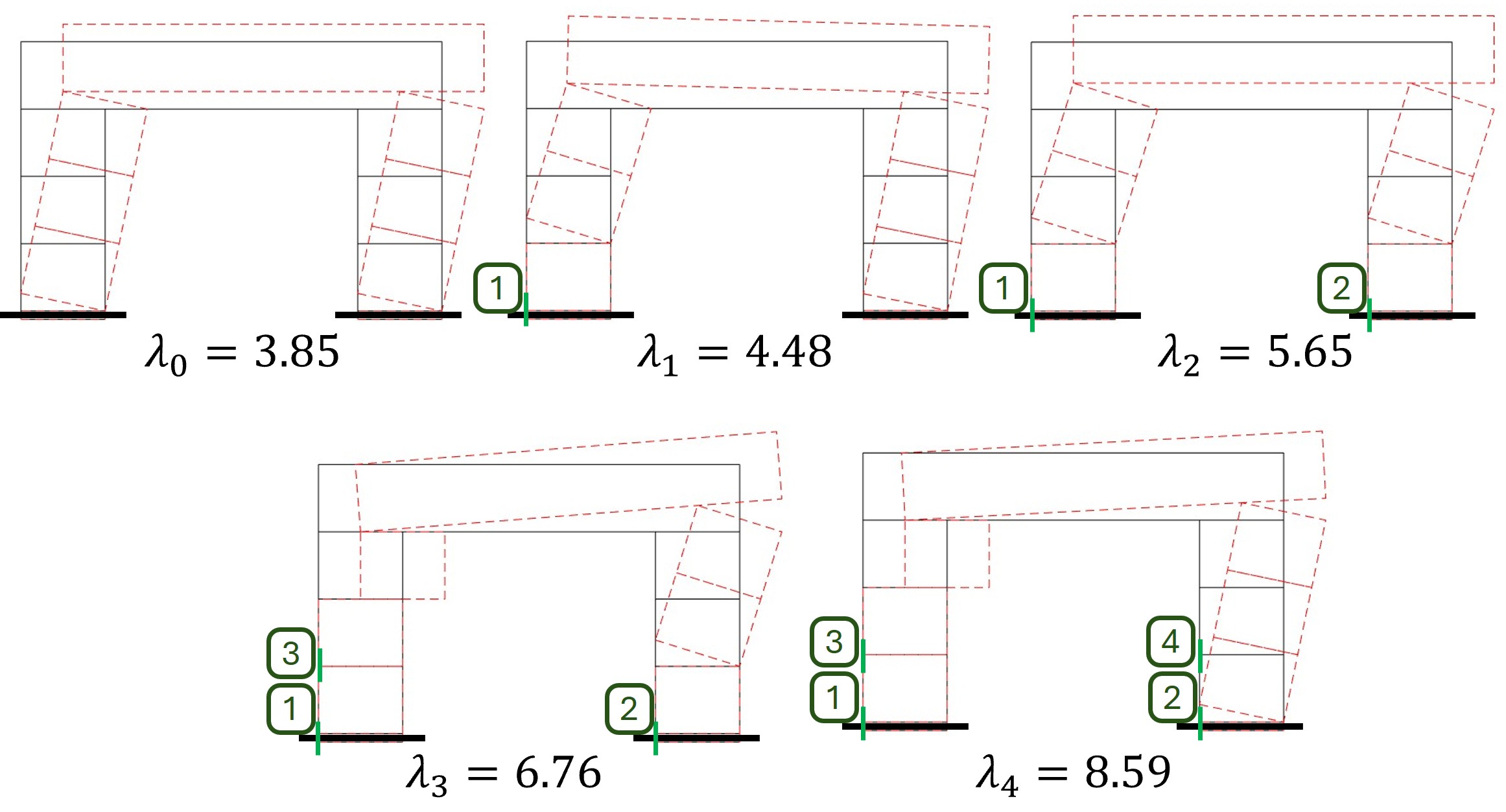}
\caption{Collapse load multiplier and associated collapse mechanism for increasing reinforcements.  \label{fig:trilite_kin}}
\end{figure}

In Figure \ref{fig:trilite_kin}, the collapse mechanisms and the corresponding collapse multipliers are presented for an increasing number of applied reinforcements. The disposition of reinforcements within the structure is also depicted in the same figure. Notice that the represented reinforced are active only for seismic inertia forces directed from left to right; since seismic actions are bidirectional, in a real-world applications also reinforcements in symmetric positions have to be applied (not depicted in Figure \ref{fig:trilite_kin}).
The first case in Figure \ref{fig:trilite_kin} represents the unreinforced structure. Subsequently, the cases obtained by adding one reinforcement at a time to the structure are represented in sequence. It is evident that the collapse multiplier increases as the number of reinforcements increases; in the last case (reinforcement applied to the first two joints starting from the bottom of both columns) the collapse load multiplier undergoes an increase of about 123\%, going from the value $\lambda = 3.85$ (unreinforced structure up) to the value $\lambda = 8.59$. Notice that the reinforced structure is still capable of forming mechanisms, but the introduction of the reinforcement generates a change in the kinematics. It is noteworthy that starting from the case where three reinforcements are applied, the kinematics does not only induce the opening of the interfaces, but also activates sliding between the blocks.

\begin{figure}[ht!]
\centering
\includegraphics[scale=0.25]{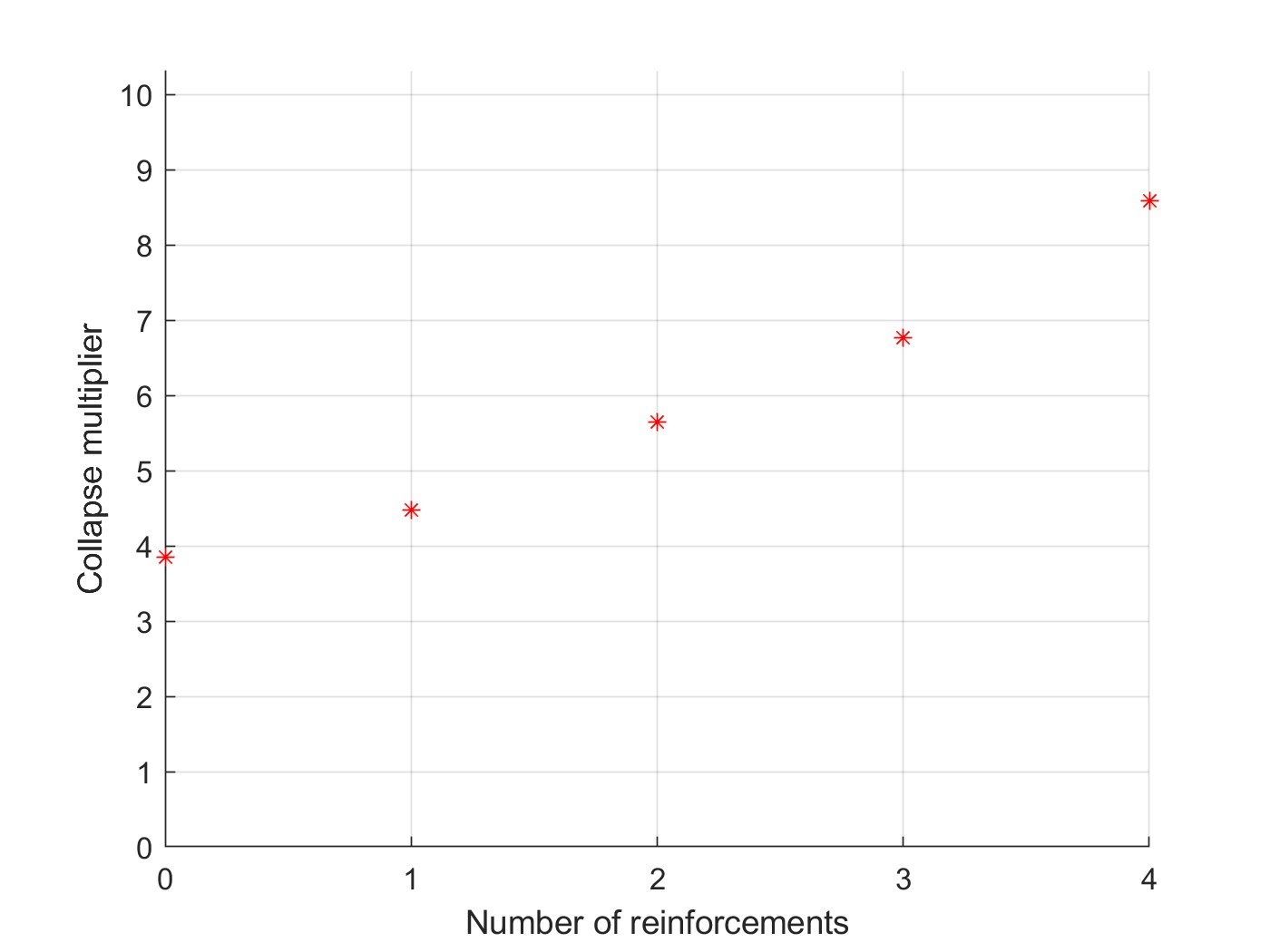}
\caption{Variation of the collapse multiplier for an increasing number of reinforcements. \label{fig:trilite_lambda}}
\end{figure}

Furthermore, it turns out that when four reinforcements are applied, an interface opening occurs at the base of the right column, despite the presence of the reinforcement. This is due to the fact that, in this case, the reinforcement reaches its limit strength and is therefore subject to plastic deformation.
Figure \ref{fig:trilite_lambda} presents the plot of the value of the collapse multiplier versus the number of applied reinforcements.

Although the application of SLA1, SLA2, and KLA leads to the same values of the collapse loads and to the same corresponding kinematics, the outcomes of these diverse limit analysis approaches yield solutions that differ from each other, particularly in the evaluation of the interface forces between blocks.

The SLA1 approach is the sole method that explicitly incorporates forces in the reinforcements as unknowns, thereby rendering it the only approach capable of evaluating such forces.
The values of these forces for different numbers of reinforcements are shown in Tab. \ref{tab:trilite_reinf}.
The evaluation of the value of forces in the reinforcements enables the initial optimization of the reinforcements' design. Indeed, the quantification of the reinforcement itself is possible based on the value of the forces arising in the reinforcements.

\begin{table}[h!]
  \begin{center}
    \caption{Forces in the reinforcements for the trilithon system.}
    \label{tab:trilite_reinf}
    \small
    \begin{tabular}{c|rrrr} 
      \textbf{Number of} &  & & & \vspace*{-2 mm}\\
      \textbf{reinforcements} & 1\;\;\;  & 2\;\;\; & 3\;\;\; & 4\;\;\;\\
      \hline
      $R_1$ & 2.1093 & 2.2677 & 5.7483 &  5.7774\\
      $R_2$ &        & 3.7194 & 3.9609 & 10.0000\\
      $R_3$ &        &        & 2.2189 &  2.0730\\
      $R_4$ &        &        &        &  3.7635
    \end{tabular}
  \end{center}
\end{table}

The data presented in Table \ref{tab:trilite_reinf} show that the yield force value is attained exclusively in reinforcement $R_2$ when four reinforcements are applied, while other the forces in other reinforcements are significantly lower than the yield force. This observation suggests that the remaining reinforcements may be designed with quite lower strength capacities (less reinforcement material quantity).

As previously mentioned, a notable benefit of masonry structures is their capacity to undergo significant differential settlement without compromising their safety and serviceability.
In general, masonry structures retain their functionality and structural safety even in the presence of cracks. This resilience is attributed to their adaptability to new configurations, akin to isostatic structures. 
It is crucial, however, to recognize that the integration of reinforcements can potentially compromise this adaptability, particularly in cases where the design is not executed in a comprehensive manner.
Therefore, a pivotal examination that must be conducted on reinforced structures pertains to the response of these structures to foundation settlements.
To investigate the effect of potential foundation settlements on the reinforced structure, an analysis is conducted to ascertain the feasibility of achieving new equilibrium configurations without inducing overstresses in the structure, under assigned support displacements. 
The KLA kinematic method is employed for this analysis, and the solution for $\lambda=0$ is derived using the total potential energy approach.

\begin{figure}[ht!]
\centering
\includegraphics[scale=0.8]{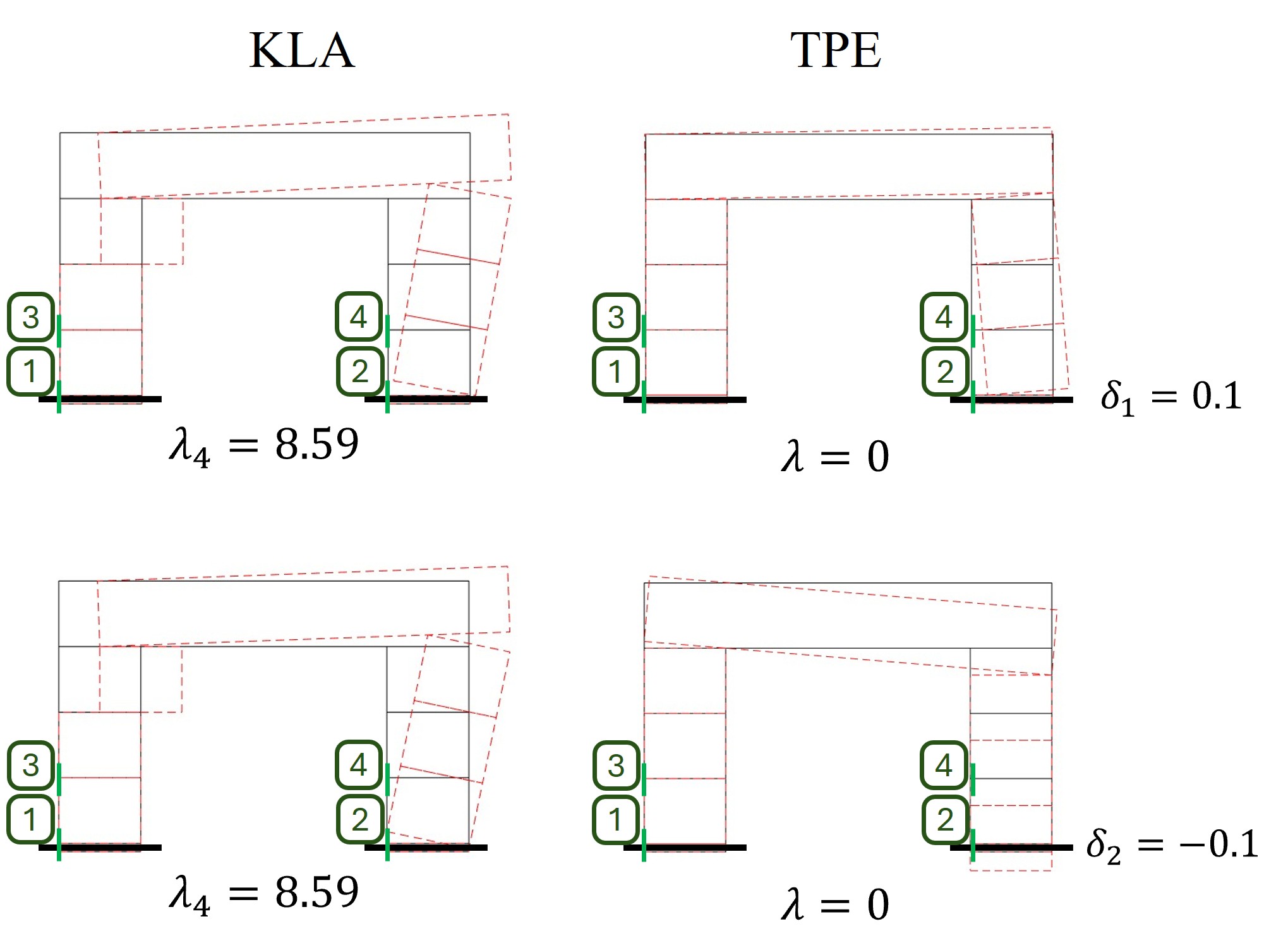}
\caption{Effect of the foundation settlement on the full reinforced structural system (displacements in the right column are much more amplified than that in the left column). \label{fig:trilite_delta}}
\end{figure}

In Figure \ref{fig:trilite_delta}, the solutions obtained by prescribing two different foundation settlements of the right pillar are reported. In particular, an assigned horizontal displacement $\delta_1=0.1$ (top row of the figure) and an assigned vertical displacement $\delta_2=-0.1$ (bottom row of the figure)  are considered for the reinforced structure. The left column of Figure \ref{fig:trilite_delta} displays the collapse load obtained by the application of the KLA approach. It is noteworthy that the collapse multiplier remains unchanged for both assigned settlements, i.e., $\delta_1=0.1$ and $\delta_2=-0.1$, when compared to the previously computed scenario with no settlement. This results is related to the assumption of small displacements. 

The right column of Figure \ref{fig:trilite_delta} displays the solution in terms of displaced configuration of the reinforced structure obtained using the Total Potential Energy (TPE) method when only dead loads are applied ($\lambda=0$). The kinematics activated by the two settlements are evident. It is noteworthy that, despite the application of reinforcements, the structure retains the capability to adapt to the foundation settlements through simple kinematic mechanisms that do not induce overstresses. Additionally, the activated kinematics is that characteristic of masonry structures, implying that the intrinsic behavior of masonry is not compromised by the reinforcements designed. In fact, the reinforcements contribute to enhancing the overall structural safety of the system (+123\% in seismic capacity) without compromising the capacity of adapting to foundation settlement and without violating the distinctive behavior of masonry structures.

\subsection{The arch}

\begin{figure}[ht!]
\centering
\includegraphics[scale=0.7]{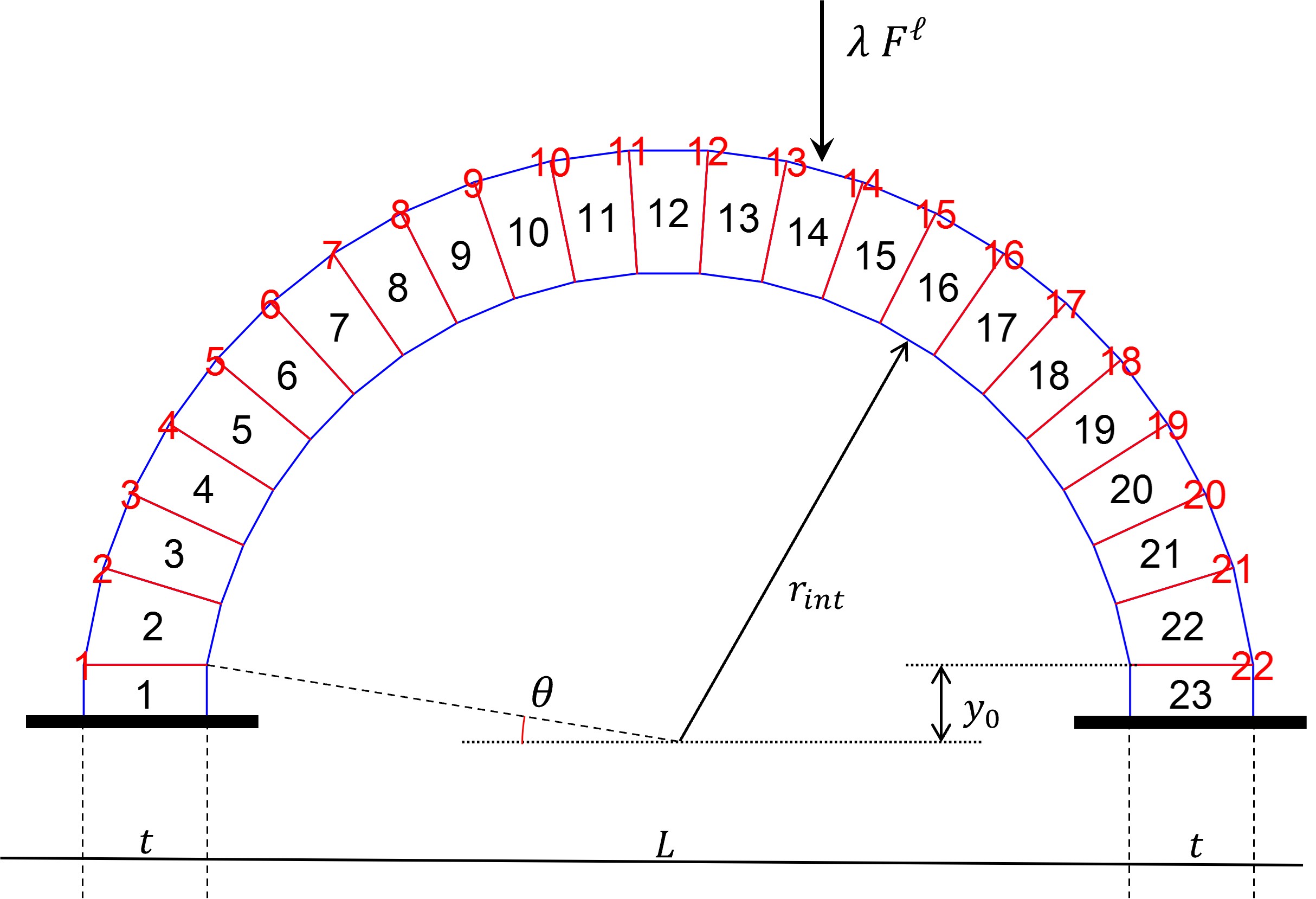}
\caption{Geometry, loading and constraints for a masonry arch. \label{fig:arch}}
\end{figure}

In this section, the response of an arch for which experimental data are available in the references \cite{Cancelliere2010,Sacco10} is investigated.
The arch is illustrated in Figure \ref{fig:arch}. 
It is an almost round arch composed of 23 clay bricks joined by mortar layers. The arch is not perfectly round, as the center of the circular mid-line, intrados and extrados lines is lower than the level of the imposts, and as the bricks acting as the imposts have horizontal faces.
The geometrical parameters that describe the arch are listed below: intrados radius $r_{int}=456$ mm, internal span $L=900$ mm, height of the bricks acting as the imposts $y_0=74$ mm, thickness of bricks $t=120$ mm, and opening angle $\theta=0.1624 \, \textrm{rad}$. The width of the arch is $w=240$ mm, and the weight of the homogenized material (i.e., brick and mortar) is $w_0=1.60e-05 \, \textrm{N/mm}^2$.
The arch is subject to a dead load resulting from its own weight and a live load consisting of a concentrated force applied at the center of block number 14, as illustrated in Figure \ref{fig:arch}.
The arch is modeled as formed by blocks joined by zero-thickness interfaces by expanding the brick dimensions until the middle of the mortar joint, and substituting the mortar with an interface of zero thickness. This interface is characterized by unilateral contact and friction, with no cohesion.
It should be noted that in Figure \ref{fig:arch}, the expanded blocks are numbered in black, while the interfaces are numbered in red.

\begin{figure}[ht!]
\centering
\includegraphics[scale=0.9]{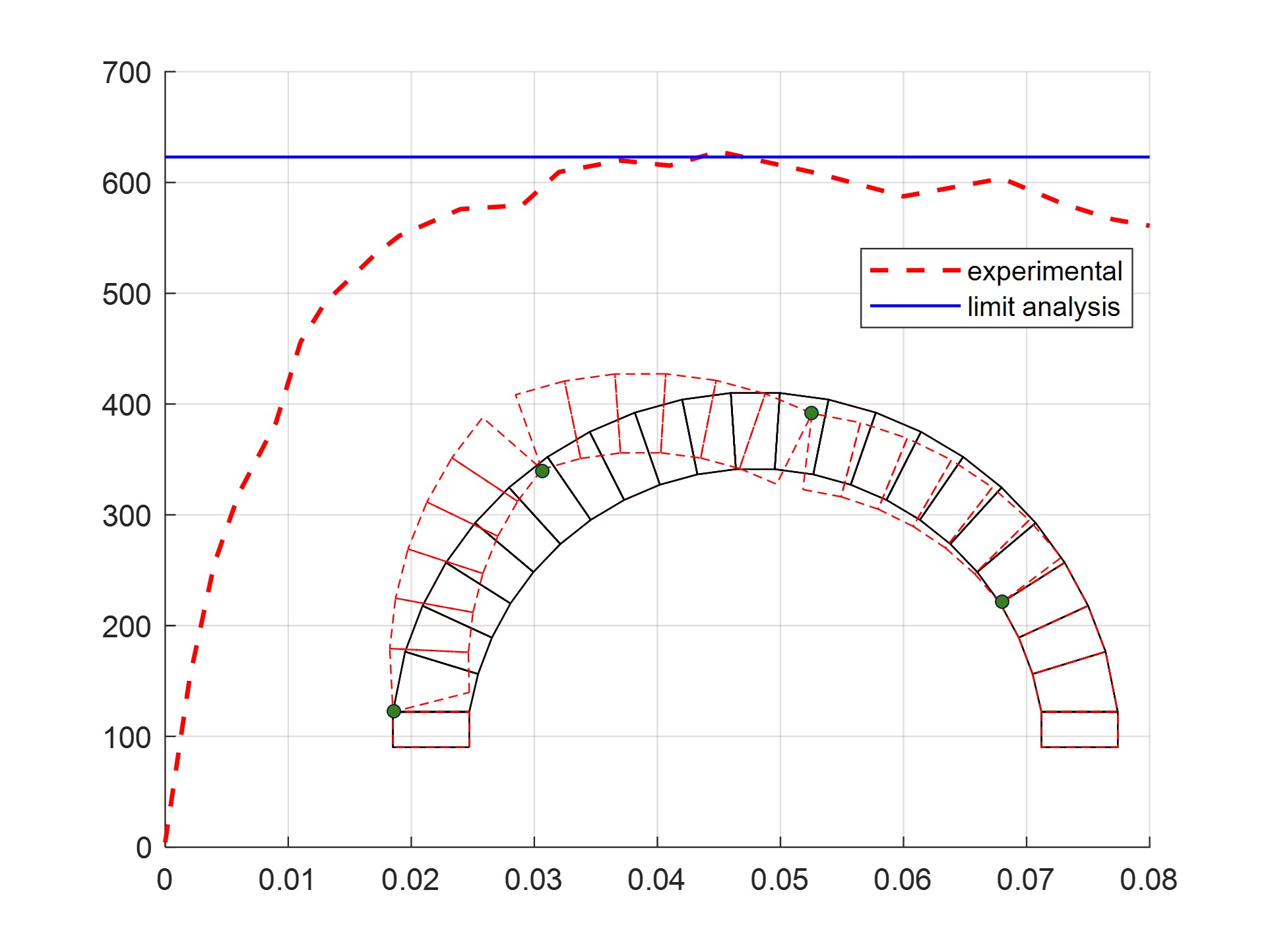}
\caption{Experimental response of the unreinforced arch and limit load evaluation. Collapse kinematics and position od the hinges. \label{fig:arch_collapse}}
\end{figure}

Initially, the limit analysis (SLA1, SLA2, or KLA approaches) is applied to determine the collapse load of the unreinforced arch. The obtained collapse load is then displayed in comparison with the response curve of the arch in terms of applied load versus vertical displacement of the loaded point, determined from the experimental test reported in \cite{Sacco10}, see Figure 18. It is evident that limit analysis provides a very accurate evaluation of the true collapse load. In the same Figure \ref{fig:arch_collapse} is depicted the collapse mechanism obtained by limit analysis, that turns out to be analogous to that observed during the experimental test (see \cite{Sacco10}, Figure 18). Also the position of the formed hinges is outlined in the figure.

Subsequently, the reinforcements are applied, and the collapse load and the traction in the reinforcements are calculated using the SLA1 approach.
The reinforcement yield strength was set at $R_y = 100 N$, a deliberately low value chosen to provide minimal, incremental strengthening of the structure.

\begin{figure}[ht!]
\centering
\includegraphics[scale=0.6]{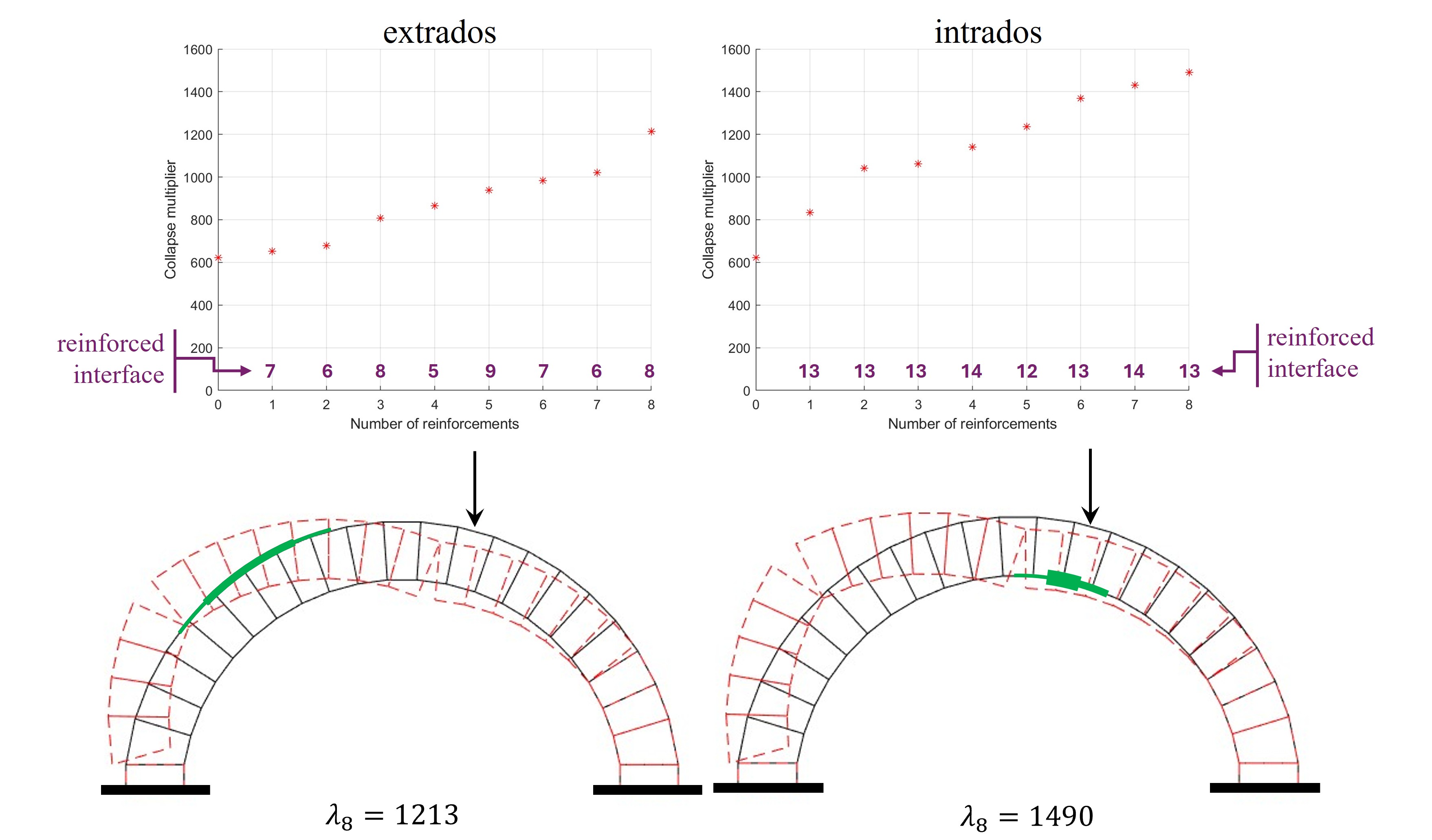}
\caption{Kinematics and collapse multiplier for the trilithon with increasing number of reinforcements. \label{fig:reinf_arch}}
\end{figure}

Figure \ref{fig:reinf_arch} presents the results for the arch with different reinforcement configurations. The left column corresponds to extrados strengthening, while the right column pertains to intrados strengthening. Each column consists of two rows: the top row displays the collapse load values as a function of the number of reinforcements applied; also the labels of the reinforced interfaces are indicated (numerical labels above the abscissa axis; the repetition of interface numbers indicates that multiple reinforcement layers were required at certain locations). The bottom row provides a schematic representation of reinforcement placement for the maximum number of reinforcements considered, and the associated collapse mechanism.

The results indicate that the optimization procedure for designing the weak reinforcement of the considered arch leads to the superposition of multiple reinforcements at the same interface, both for extrados and intrados applications. Indeed, for extrados strengthening the sequence of the reinforced interfaces is 7, 6, 8, 5, 9, again 7, again 6, again 8. For for intrados strengthening the sequence of the reinforced interfaces is three times 13, 14, 12, again 13, again 14, again 13.
This is evident in the maximum reinforcement placement schematics depicted in the figure, where certain reinforcements exhibit higher thickness.
The collapse load graphs for both extrados and intrados strengthening show that in the first case, by reinforcing only interfaces from 5 to 9 is possible about doubling (from 623 to 1213) the collapse load, while for intrados strengthening by reinforcing only interfaces from 12 to 14 the collapse load increases about 2.4 times  (from 623 to 1490). This, without compromising the free hinging capacity of not reinforced interfaces, and therefore the adaptability capacity of the arch to foundation settlements and its distinctive behavior as a masonry structure. In both cases the number of reinforcements (8) is the same.

However, it is well-known that intrados reinforcements suffer from significant debonding issues from the substrate, making them generally less preferred with respect to extrados strengthening. Indeed, extrados reinforcements — due to geometric reasons — tend to remain adherent to the arch when subjected to tensile forces. Moreover, extrados strengthening is often preferred since the historical constructions the intrados the texture of the masonry must be left visible on the intrados and/or there are paintings on the intrados that must be preserved.

Tables \ref{tab:arch_reinf_ex} and \ref{tab:arch_reinf_in} present for the reinforcements placed at extrados and at intrados, respectively, the force distribution within the reinforcements as a function of their quantity. Notably, the limiting force $R_y=100.000 N$ is attained in only one reinforcement for both extrados and intrados reinforcement systems when 8 reinforcements are applied.
This suggests that, for an optimal design, the number of reinforcements could be reduced in regions where the force remains below $R_y$  without necessarily diminishing the critical load capacity.

\begin{table}[h!]
  \begin{center}
    \caption{Forces in the reinforcements for the arch reinforced at the extrados.}
    \label{tab:arch_reinf_ex}
    \small
    \begin{tabular}{c|rrrrrrrr} 
      \textbf{Number of} &  & & & \vspace*{-2 mm}\\
      \textbf{reinforcements} & 1\;\;\;  & 2\;\;\; & 3\;\;\; & 4\;\;\; & 5\;\;\; & 6\;\;\; & 7\;\;\; & 8\;\;\;\\
      \hline
      $R_1$ & 9.545 &  17.493 & 58.896 & 77.019 & 100.000 &  57.237 &  63.148 &  93.605 \\
      $R_2$ &        &  7.5748 & 47.038 & 64.311 &  86.215 & 100.000 &  55.643 &  84.669\\
      $R_3$ &        &         & 38.054 & 54.710 &  75.832 &  89.125 & 100.000 &  77.993\\
      $R_4$ &        &         &        & 14.124 &  32.034 &  43.305 &  52.530 & 100.000\\
      $R_5$ &        &         &        &        &  16.360 &  26.656 &  35.078 &  78.441\\
      $R_6$ &        &         &        &        &         &  57.237 &  63.148 &  93.605\\
      $R_7$ &        &         &        &        &         &         &  55.643 &  84.669\\
      $R_8$ &        &         &        &        &         &         &         &  77.993
    \end{tabular}
  \end{center}
\end{table}

\begin{table}[h!]
  \begin{center}
    \caption{Forces in the reinforcements for the arch reinforced at the intrados.}
    \label{tab:arch_reinf_in}
    \small
    \begin{tabular}{c|rrrrrrrr} 
      \textbf{Number of} &  & & & \vspace*{-2 mm}\\
      \textbf{reinforcements} & 1\;\;\;  & 2\;\;\; & 3\;\;\; & 4\;\;\; & 5\;\;\; & 6\;\;\; & 7\;\;\; & 8\;\;\;\\
      \hline
      $R_1$ & 100.00 & 100.00 & 70.104 & 83.450 & 100.000 &  91.921 & 100.000 &  86.085 \\
      $R_2$ &        & 100.00 & 70.104 & 83.450 & 100.000 &  91.921 & 100.000 &  86.085\\
      $R_3$ &        &        & 70.104 & 83.450 & 100.000 &  91.921 & 100.000 &  86.085\\
      $R_4$ &        &        &        & 25.442 &  56.992 & 100.000 &  60.268 &  69.933\\
      $R_5$ &        &        &        &        &  27.572 &  65.160 &  83.105 & 100.000\\
      $R_6$ &        &        &        &        &         &  91.921 & 100.000 &  86.085\\
      $R_7$ &        &        &        &        &         &         &  60.268 &  69.993\\
      $R_8$ &        &        &        &        &         &         &         &  86.085
    \end{tabular}
  \end{center}
\end{table}

The bottom panel of Figure \ref{fig:reinf_arch} shows the collapse mechanism for both configurations (extrados and intrados strengthening), whether the 8 reinforcements are applied to the extrados or intrados, revealing a very similar failure behavior in both cases.
Looking at that figure, it can be noted that for extrados reinforcement systems, the collapse mechanism develops the opening of a hinge at interface No. 4, which is a reinforced interface. As demonstrated in Table \ref{tab:arch_reinf_ex}, reinforcement No. 4 reaches its limiting force capacity, allowing interface opening to occur.
Similarly, for intrados reinforcement systems, Table \ref{tab:arch_reinf_in} shows that the limiting force is attained at reinforcement No. 5, which opens in the collapse mechanism.

As can also be observed in Figure \ref{fig:reinf_arch}, the weak reinforcement design yields a very limited strengthening intervention, counteracting the opening of only a few interfaces. It is therefore expected that the structure's ability to adapt to potential foundation settlements remains unchanged. For assessing this ability, as in the previous example, the limit load computation is performed in the presence of foundation settlements, and, moreover, the structure subjected solely to vertical loads and settlements is analyzed through the TPE approach. Notice that, analogously to the trilite case, the collapse load turns out to be uninfluenced by the foundation settlements imposed.

\begin{figure}[ht!]
\centering
\includegraphics[scale=0.6]{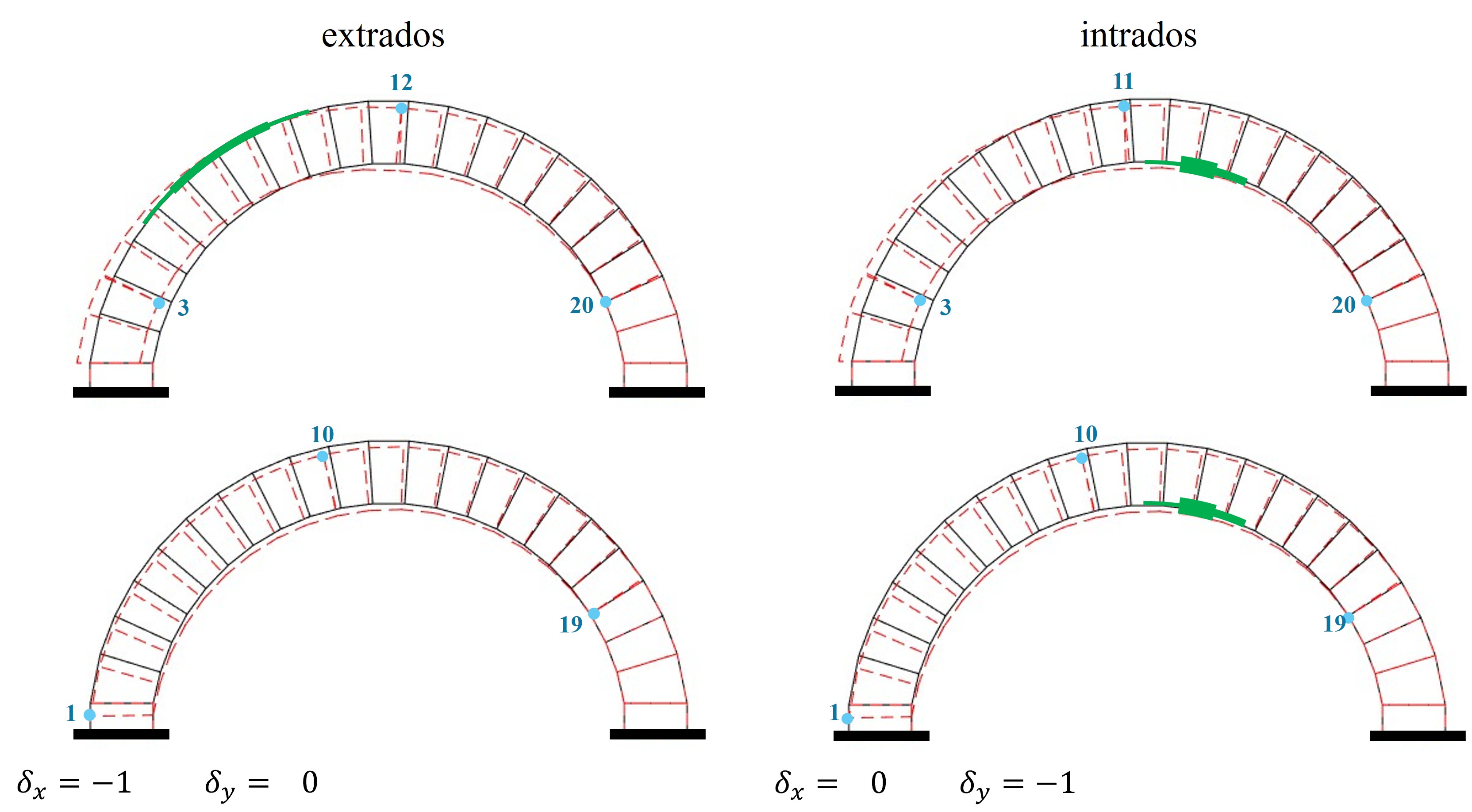}
\caption{Kinematic mechanisms for the fully reinforced arch (with reinforcement applied at either extrados or intrados) under combined loading conditions: self-weight distribution and prescribed settlement at the right springing of the masonry arch. \label{fig:arch_TPE}}
\end{figure}

Figure \ref{fig:arch_TPE} illustrates the mechanisms of adaptation to foundation settlements of the arch reinforced at either the extrados (left) or intrados (right). In particular, calculations have been performed considering the arch subjected solely to self-weight loading and to a prescribed displacement of the left springing in either horizontal direction ($\delta_x=-1$, upper panel of the figure) or vertical direction ($\delta_y=-1$, lower panel of the figure).

The resulting mechanisms exhibit relatively small displacements (amplified 25$\times$ in the figure for making them visible), with the displaced configuration remaining notably close to the undeformed state. Through the mechanism analysis and numerical results, hinge locations are identified (marked by blue circles in the figure, labeled with their corresponding interface numbers).
In all cases, three hinges form that, combined with the prescribed foundation settlement, trigger an isostatic behavior of the arch and allow to form a mechanism that avoid the emergence of additional stresses in the structure, thus preserving the masonry arch's key adaptive capacity to foundation settlements without inducing overstressing.

\section{Conclusions}

The proposed approach aims at the design of minimally invasive, strategically placed reinforcements that enhance structural capacity without compromising the adaptive and isostatic characteristics behavior of masonry structures. This approach is called "weak reinforcement" since it consists in the selective reinforcement of only some suitable part of the structure, leaving the remaining part free to behave as conventional masonry structures, where the small tensile strength involve unilateral behavior of joints, capable of forming hinges and also, for special conditions, of sliding.
The concept of weak reinforcement is somewhat opposed as that of conventional strengthening interventions, where the reinforcement is applied to the whole face of a masonry structural element, or eventually to both the whole faces of the element. The difference does not concern only the extension of the reinforcement, but the function thought for the reinforcement itself. Indeed, the conventional way of thinking is that the more compensation is provided for the limited tensile strength of the masonry by preventing the opening of hinges at joints through the application of reinforcing materials, the more effective the reinforcement will be. This approach might appear effective based on the results of some simplified calculations, but such way of strengthening masonry construction, if can improve the structural capacity with respect to some special class of action, on the other side can weaken, even significantly, the structure. 
Indeed, since hinge openings are totally prevented, the structure cannot undergo rocking motions, they can make the structure survive even very strong earthquakes. Moreover, the structure looses the ability of adapting to support displacements and settlements typical of masonry construction: the conventionally strengthened structure cannot behave isostatically, and therefore even small support displacements can produce high stress states in the masonry. Finally, a conventionally strengthened structural element generally transmits very high actions to the surrounding, not reinforced structural elements, triggering weak-ring-of-the-chain collapse mechanisms.
Another aspect to be taken into account is that often with conventional strengthening techniques is not possible to graduate the reinforcement to obtain the desired capacity level, and it is necessary to accept a capacity higher than that required. 
All the above considerations have lead to the search of a strengthening approach capable of leaving to the masonry structure the capacity of forming mechanisms and of adapting isostatically to support displacements and settlements. Moreover, this study has been performed also with the objective of optimizing the reinforcement quantity for not exceeding the desired capacity, therefore limiting the transmission of undesired forces to the not reinforced surrounding parts or the structures, and for perturbing as little as possible the construction and its typical masonry behavior.
In short, the proposed new type of "weak reinforcement" strategy aims at selectively and strategically introducing in only some joints reinforcements acting like "reinforcing forces", limited by the yield strength of the material, rather than perfect kinematic constraints. This approach appears to be in line with the criteria in \cite{ICOMOS2024} for the structural restoration of architectural heritage: indeed, along with the minimum intervention, the objective of the structural improvement (performance level suitably lesser than that of a new structure) and of the preservation of the original structural behavior is pursued.

To achieve this goal, and to optimize the reinforcement quantity, it was necessary to develop a specific structural analysis procedure based on static and kinematic non-standard limit analysis, capable of taking into account also the behavior of non-associative joints, and the presence of reinforcements preventing the separation of joints ends. Moreover, based on a suitable approach for the minimization of the total potential energy, the ability of the structure to adapt to support displacement can be studied. 

The integration of limit analysis, numerical optimization, and energy minimization provides a robust tool for practical engineering applications, especially in the conservation of historic masonry.

Notice that this paper does not aim to study how such a strengthening system can be realized from a technological point of view, but rather to define a new strengthening strategy, having objectives very different from the conventional one, and aimed at the conservation of masonry buildings not only by minimizing the amount of strengthening material used, but above all by the fact that the increase in capacity is obtained by preserving the original "masonry-like" behavior.
The obtained results encourage to continue further in this direction, and pave the way for future developments.

First, the proposed weak reinforcement approach has to be studied from the point of view of the technological feasibility, for example suitably using FRP strips with mechanical anchors fixed astride the joints. Moreover, numerical results have to be experimentally validated by laboratory tests.
Another significant investigation direction concerns the study of "elastic" weak reinforcements, capable of exerting forces at the joint ends that are proportional to joint openings, and therefore allowing for joint hinging from the beginning, while supplying tensile forces beneficial for structural capacity.
Finally, among the wide panorama of possible research directions for developing the weak reinforcement concept, a crucial issue to be investigated is that of structures where three-dimensional effects are non negligible. Indeed, the theoretical framework here proposed is conceived only for 2D structures, and the extension to 3D problems might be not trivial.

\section*{Acknowledgment}

\bibliographystyle{ieeetr}
\bibliography{biblio}

\begin{thebibliography}{10}

\bibitem{Drysdale1999}
R.~G. Drysdale, A.~A. Hamid, and L.~R. Baker, {\em Masonry Structures: Behavior and Design}.
\newblock The Masonry Society, 2~ed., 1999.

\bibitem{Hendry2004}
A.~W. Hendry, B.~P. Sinha, and S.~R. Davies, {\em Design of Masonry Structures}.
\newblock Taylor \& Francis, 3~ed., 2004.

\bibitem{Como2016}
M.~Como, {\em Statics of Historic Masonry Constructions}.
\newblock Springer International Publishing, 2016.

\bibitem{Heyman1997}
J.~Heyman, {\em The Stone Skeleton: Structural Engineering of Masonry Architecture}.
\newblock Cambridge University Press, 1997.

\bibitem{Portioli2016}
F.~Portioli and L.~Cascini, ``Assessment of masonry structures subjected to foundation settlements using rigid block limit analysis,'' {\em Engineering Structures}, vol.~113, pp.~347 -- 361, 2016.

\bibitem{Tralli2020}
A.~Tralli, A.~Chiozzi, N.~Grillanda, and G.~Milani, ``Masonry structures in the presence of foundation settlements and unilateral contact problems,'' {\em International Journal of Solids and Structures}, vol.~191-192, pp.~187--201, 2020.

\bibitem{Iannuzzo2025}
A.~Iannuzzo and V.~Mallardo, ``A novel approach to model differential settlements and crack patterns in masonry structures,'' {\em Engineering Structures}, vol.~323, no.~Part A, p.~119220, 2025.

\bibitem{Ferretti2006}
D.~Ferretti and Z.~P. Baz˘ant, ``Stability of ancient masonry towers: Stress redistribution due to drying, carbonation, and creep,'' {\em Cement and Concrete Research}, vol.~36, no.~7, pp.~1389--1398, 2006.

\bibitem{Cecchi2012}
A.~Cecchi and A.~Tralli, ``A homogenized viscoelastic model for masonry structures,'' {\em International Journal of Solids and Structures}, vol.~49, no.~13, pp.~1485--1496, 2012.

\bibitem{Sánchez-Beitia2017}
S.~Sánchez-Beitia, D.~Luengas-Carreño, and M.~Crespo~de Antonio, ``The presence of secondary creep in historic masonry constructions: A hidden problem,'' {\em Engineering Failure Analysis}, vol.~82, pp.~315--326, 2017.

\bibitem{Tomazevic1999}
M.~Tomazevič, {\em Earthquake-Resistant Design of Masonry Buildings}.
\newblock Imperial College Press, July 1999.

\bibitem{Asteris2014}
P.~G. Asteris, M.~P. Chronopoulos, C.~Z. Chrysostomou, H.~Varum, V.~Plevris, N.~Kyriakides, and V.~Silva, ``Seismic vulnerability assessment of historical masonry structural systems,'' {\em Engineering Structures}, vol.~62–63, pp.~118--134, 2014.

\bibitem{Lagomarsino2021}
S.~Lagomarsino, S.~Cattari, and D.~Ottonelli, ``The heuristic vulnerability model: fragility curves for masonry buildings,'' {\em Bulletin of Earthquake Engineering}, vol.~19, pp.~3129--–3163, 2021.

\bibitem{ICOMOS2024}
I.~. ISCARSAH, {\em Guidelines on the Analysis, Conservation and Structural Restoration of Architectural Eeritage}.
\newblock 2024.

\bibitem{Grande2008}
E.~Grande, G.~Milani, and E.~Sacco, ``Modelling and analysis of frp-strengthened masonry panels,'' {\em Engineering Structures}, vol.~30, no.~7, pp.~1842 -- 1860, 2008.

\bibitem{Cancelliere2010}
I.~Cancelliere, M.~Imbimbo, and E.~Sacco, ``Experimental tests and numerical modeling of reinforced masonry arches,'' {\em Engineering Structures}, vol.~32, no.~3, pp.~776 -- 792, 2010.

\bibitem{Carozzi2018}
F.~G. Carozzi, C.~Poggi, E.~Bertolesi, and G.~Milani, ``Ancient masonry arches and vaults strengthened with trm, srg and frp composites: Experimental evaluation,'' {\em Composite Structures}, vol.~187, pp.~466 -- 480, 2018.

\bibitem{Alecci2016}
V.~Alecci, G.~Misseri, L.~Rovero, G.~Stipo, M.~De~Stefano, L.~Feo, and R.~Luciano, ``Experimental investigation on masonry arches strengthened with pbo-frcm composite,'' {\em Composites Part B: Engineering}, vol.~100, pp.~228 -- 239, 2016.

\bibitem{Leone2017}
M.~Leone, M.~Aiello, A.~Balsamo, F.~Carozzi, F.~Ceroni, M.~Corradi, M.~Gams, E.~Garbin, N.~Gattesco, P.~Krajewski, C.~Mazzotti, D.~Oliveira, C.~Papanicolaou, G.~Ranocchiai, F.~Roscini, and D.~Saenger, ``Glass fabric reinforced cementitious matrix: Tensile properties and bond performance on masonry substrate,'' {\em Composites Part B: Engineering}, vol.~127, pp.~196 -- 214, 2017.

\bibitem{Meriggi2022}
P.~Meriggi, C.~Caggegi, A.~Gabor, and G.~de~Felice, ``Shear-compression tests on stone masonry walls strengthened with basalt textile reinforced mortar (trm),'' {\em Construction and Building Materials}, vol.~316, p.~125804, 2022.

\bibitem{Castellano2025}
A.~Castellano, A.~Fraddosio, F.~Paparella, M.~Piccioni, and T.~Kundu, ``The evaluation of the adhesion defects in frcm reinforcements for masonry constructions by sideband peak count based nonlinear acoustic technique,'' {\em JVC/Journal of Vibration and Control}, vol.~31, no.~3-4, pp.~258 -- 270, 2025.

\bibitem{Cimellaro2011}
G.~Cimellaro, A.~Reinhorn, and A.~De~Stefano, ``Introspection on improper seismic retrofit of basilica santa maria di collemaggio after 2009 italian earthquake,'' {\em International Journal of Architectural Heritage}, vol.~10, p.~153–161, 2011.

\bibitem{Arcidiacono2016}
V.~Arcidiacono, G.~P. Cimellaro, E.~Piermarini, and J.~Ochsendorf, ``The dynamic behavior of the basilica of san francesco in assisi using simplified analytical models,'' {\em International Journal of Architectural Heritage}, vol.~10, no.~7, p.~938–953, 2016.

\bibitem{Borri2019}
A.~Borri and M.~Corradi, ``Architectural heritage: A discussion on conservation and safety,'' {\em Heritage}, vol.~2, no.~1, pp.~631--647, 2019.

\bibitem{Grande2011}
E.~Grande, M.~Imbimbo, and E.~Sacco, ``A beam finite element for nonlinear analysis of masonry elements with or without fiber-reinforced plastic (frp) reinforcements,'' {\em International Journal of Architectural Heritage}, vol.~5, no.~6, pp.~693 -- 716, 2011.

\bibitem{Grande2023}
E.~Grande, M.~Gafone, T.~Rotunno, and G.~Milani, ``Modeling of shear-lap tests of flat and curved masonry specimens strengthened by frcm,'' {\em Structures}, vol.~52, pp.~437 -- 448, 2023.

\bibitem{Castellano2023}
A.~Castellano, A.~Fraddosio, D.~Oliveira, M.~Piccioni, E.~Ricci, and E.~Sacco, ``The evaluation of the adhesion defects in frcm reinforcements for masonry constructions by sideband peak count based nonlinear acoustic technique,'' {\em JVC/Journal of Vibration and Control}, vol.~31, no.~3-4, pp.~258 -- 270, 2025.

\bibitem{Block2006}
P.~Block, T.~Ciblac, and J.~Ochsendorf, ``Real-time limit analysis of vaulted masonry buildings,'' {\em Computers and Structures}, vol.~84, no.~29-30, pp.~1841 -- 1852, 2006.

\bibitem{Milani2006}
G.~Milani, P.~Lourenço, and A.~Tralli, ``Homogenised limit analysis of masonry walls, part ii: Structural examples,'' {\em Computers and Structures}, vol.~84, no.~3-4, pp.~181 -- 195, 2006.

\bibitem{Brandonisio2020}
G.~Brandonisio, M.~Angelillo, and A.~De~Luca, ``Seismic capacity of buttressed masonry arches,'' {\em Engineering Structures}, vol.~215, p.~110661, 2020.

\bibitem{Fraddosio2020}
A.~Fraddosio, N.~Lepore, and M.~D. Piccioni, ``Thrust surface method: An innovative approach for the three-dimensional lower bound limit analysis of masonry vaults,'' {\em Engineering Structures}, vol.~202, p.~109846, 2020.

\bibitem{Livesley1978}
R.~Livesley, ``Limit analysis of structures formed from rigid blocks,'' {\em International Journal for Numerical Methods in Engineering}, vol.~12, no.~12, pp.~1853 -- 1871, 1978.

\bibitem{Alexakis2015}
H.~Alexakis and N.~Makris, ``Limit equilibrium analysis of masonry arches,'' {\em Archive of Applied Mechanics}, vol.~85, no.~9-10, pp.~1363 -- 1381, 2015.

\bibitem{Chiozzi2017}
A.~Chiozzi, G.~Milani, and A.~Tralli, ``Fast kinematic limit analysis of frp-reinforced masonry vaults. ii: Numerical simulations,'' {\em Journal of Engineering Mechanics}, vol.~143, no.~9, p.~04017072, 2017.

\bibitem{Michiels2017}
T.~Michiels, R.~Napolitano, S.~Adriaenssens, and B.~Glisic, ``Comparison of thrust line analysis, limit state analysis and distinct element modeling to predict the collapse load and collapse mechanism of a rammed earth arch,'' {\em Engineering Structures}, vol.~148, pp.~145 -- 156, 2017.

\bibitem{Fortunato2018}
o.~A. Fortunat, F.~Fabbrocino, M.~Angelillo, and F.~Fraternali, ``Limit analysis of masonry structures with free discontinuities,'' {\em Meccanica}, vol.~53, no.~7, pp.~1793 -- 1802, 2018.

\bibitem{Iannuzzo2021}
A.~Iannuzzo, A.~Dell'Endice, T.~Van~Mele, and P.~Block, ``Numerical limit analysis-based modelling of masonry structures subjected to large displacements,'' {\em Computers and Structures}, vol.~242, p.~106372, 2021.

\bibitem{Nodargi2022}
A.~Nodargi and P.~Bisegna, ``Generalized thrust network analysis for the safety assessment of vaulted masonry structures,'' {\em Engineering Structures}, vol.~270, p.~114878, 2022.

\bibitem{Nodargi2023}
A.~Nodargi, ``An isogeometric collocation method for the static limit analysis of masonry domes under their self-weight,'' {\em Computer Methods in Applied Mechanics and Engineering}, vol.~416, p.~116375, 2023.

\bibitem{Delpiero1998}
G.~Del~Piero, ``Limit analysis and no-tension materials,'' {\em International Journal of Plasticity}, vol.~14, no.~1-3, pp.~259--271, 1998.

\bibitem{Sacchi1968}
G.~Sacchi and M.~Save, ``A note on the limit loads of non-standard materials,'' {\em Meccanica}, vol.~3, pp.~43--45, 1968.

\bibitem{Radenkovic1961}
D.~Radenkovic, {\em Théorie des Charges Limitées, Extension à la Mécanique des Sols}.
\newblock “Seminaires de Plasticité Ec. Polytech.”, Pub. Sc. et Tech. Minist. Air. n. N. T. 116, 1961.

\bibitem{Salençon1973}
J.~Salençon, ``Théorèmes généraux de l'analyse limite,'' {\em Revue Française de Mécanique}, vol.~47, pp.~9--15, 1973.

\bibitem{DeSaxcé1998}
G.~de~Saxcé and L.~Bousshine, ``Limit analysis theorems for implicit standard materials: Application to the unilateral contact with dry friction and the non-associated flow rules in soils and rocks,'' {\em International Journal of Mechanical Sciences}, vol.~40, pp.~387--398, 1998.

\bibitem{Gilbert2006}
M.~Gilbert, C.~Casapulla, and H.~M. Ahmed, ``Limit analysis of masonry block structures with non-associative frictional joints using linear programming,'' {\em Computers and Structures}, vol.~84, no.~13--14, pp.~873--887, 2006.

\bibitem{Trentadue2013}
F.~Trentadue and G.~Quaranta, ``Limit analysis of frictional block assemblies by means of fictitious associative-type contact interface laws,'' {\em International Journal of Mechanical Sciences}, vol.~70, pp.~140--145, 2013.

\bibitem{Portioli2014}
F.~Portioli, C.~Casapulla, M.~Gilbert, and L.~Cascini, ``Limit analysis of 3d masonry block structures with non-associative frictional joints using cone programming,'' {\em Computers and Structures}, vol.~143, pp.~108 -- 121, 2014.

\bibitem{Portioli2015}
F.~Portioli, C.~Casapulla, and L.~Cascini, ``An efficient solution procedure for crushing failure in 3d limit analysis of masonry block structures with non-associative frictional joints,'' {\em International Journal of Solids and Structures}, vol.~69-70, pp.~252 -- 266, 2015.

\bibitem{Nodargi2019AVF}
N.~Nodargi, C.~Intrigila, and P.~Bisegna, ``A variational-based fixed-point algorithm for the limit analysis of dry-masonry block structures with non-associative coulomb friction,'' {\em International Journal of Mechanical Sciences}, 2019.

\bibitem{Mousavian2020}
E.~Mousavian and C.~Casapulla, ``The role of different sliding resistances in limit analysis of hemispherical masonry domes,'' {\em Frattura ed Integrita Strutturale}, vol.~14, no.~51, pp.~336 -- 355, 2020.

\bibitem{Hua2022}
Y.~Hua and G.~Milani, ``Rigid block limit analysis of masonry arches with associated and non-associated slides,'' in {\em From Corbel Arches to Double Curvature Vaults Analysis}, pp.~169–--203, Springer Cham, 2022.

\bibitem{Cocchetti2024}
G.~Cocchetti and E.~Rizzi, ``Finite-friction least-thickness self-standing domains of symmetric circular masonry arches,'' {\em Structures}, vol.~66, p.~106800, 2024.

\bibitem{Gagliardo2024}
R.~Gagliardo, L.~Cascini, and F.~Portioli, ``3d rigid block limit analysis of a masonry cross vault subjected to shaking table tests,'' {\em International Journal of Architectural Heritage}, vol.~18, no.~12, pp.~1886 -- 1899, 2024.

\bibitem{Rios2025}
A.~Rios, B.~Nela, M.~Pingaro, E.~Reccia, and P.~Trovalusci, ``Parametric analysis of masonry arches following a limit analysis approach: Influence of joint friction, pier texture, and arch shallowness,'' {\em Mathematics and Mechanics of Solids}, vol.~30, no.~1, pp.~137 -- 165, 2025.

\bibitem{trentadue2025}
F.~Trentadue, D.~De~Tommasi, and N.~Marasciuolo, ``A new approach to the limit analysis of masonry structures as assemblies of rigid-plastic blocks with frictional sliding contacts,'' {\em Structures}, vol.~76, p.~108895, 2025.

\bibitem{Nodargi2021}
A.~Nodargi and P.~Bisegna, ``Collapse capacity of masonry domes under horizontal loads: A static limit analysis approach,'' {\em International Journal of Mechanical Sciences}, vol.~212, p.~106827, 2021.

\bibitem{Caporale2006}
A.~Caporale, R.~Luciano, and L.~Rosati, ``Limit analysis of masonry arches with externally bonded frp reinforcements,'' {\em Computer Methods in Applied Mechanics and Engineering}, vol.~196, no.~1-3, pp.~247 -- 260, 2006.

\bibitem{Fabbrocino2015}
F.~Fabbrocino, I.~Farina, V.~Berardi, A.~Ferreira, and F.~Fraternali, ``On the thrust surface of unreinforced and frp-/frcm-reinforced masonry domes,'' {\em Composites Part B: Engineering}, vol.~83, pp.~297 -- 305, 2015.

\bibitem{Nela2022}
B.~Nela, A.~Jiménez~Rios, M.~Pingaro, E.~Reccia, and P.~Trovalusci, ``Limit analysis of locally reinforced masonry arches,'' {\em Engineering Structures}, vol.~271, p.~114921, 2022.

\bibitem{Sacco10}
E.~Sacco and J.~Toti, ``Interface elements for the analysis of masonry structures,'' {\em International Journal for Computational Methods in Engineering Science and Mechanics}, vol.~11, no.~6, pp.~354--373, 2010.

\bibitem{sloan1995}
S.~W. Sloan and P.~W. Kleeman, ``Upper bound limit analysis using discontinuous velocity fields,'' {\em Computer Methods in Applied Mechanics and Engineering}, vol.~127 (1–4), pp.~293--314, Nov. 1995.

\bibitem{Grillanda2025}
N.~Grillanda and V.~Mallardo, ``Compatible strain-based upper bound limit analysis model for masonry walls under in-plane loading,'' {\em Computers and Structures}, vol.~313, p.~107743, 2025.

\end{thebibliography}

\end{document}